\documentstyle{amsart}
\makeatletter
\newskip\dgARROWLENGTH  \dgARROWLENGTH=2.5em\relax
\newskip\dgHORIZPAD     \dgHORIZPAD=1em\relax
\newskip\dgVERTPAD      \dgVERTPAD=2ex\relax
\newskip\dgLABELOFFSET  \dgLABELOFFSET=.7ex\relax
\newcount\dgARROWPARTS  \dgARROWPARTS=4\relax
\newcommand{\dgeverynode}{\displaystyle}
\newcommand{\dgeverylabel}{\scriptstyle}
\newskip\dgDOTSPACING   \dgDOTSPACING=0.35em
\newskip\dgDOTSIZE      \dgDOTSIZE=1.5\fontdimen8\tenln
\newcount\dgMAXSQUARE   \dgMAXSQUARE=4\relax
\newcount\dgMINDBLSQ    \dgMINDBLSQ=10\relax
\newskip\dgCOLUMNWIDTH  \dgCOLUMNWIDTH=2em\relax
\chardef\f@ur=4
\@namedef{dgo@}{\let\dg@VECTOR=\vector
   \dg@LBLPOS=\dgARROWPARTS \divide\dg@LBLPOS\tw@}%
\@namedef{dgo@-}{\let\dg@VECTOR=\line}%
\@namedef{dgo@!}{\let\dg@VECTOR=\dg@novector}%
\@namedef{dgo@1}{\dg@LBLPOS=\@ne\relax}%
\@namedef{dgo@2}{\dg@LBLPOS=\tw@\relax}%
\@namedef{dgo@3}{\dg@LBLPOS=\thr@@\relax}%
\@namedef{dgo@4}{\dg@LBLPOS=\f@ur\relax}%
\@namedef{dgo@5}{\dg@LBLPOS=5\relax}%
\@namedef{dgo@6}{\dg@LBLPOS=6\relax}%
\@namedef{dgo@7}{\dg@LBLPOS=7\relax}%
\@namedef{dgo@8}{\dg@LBLPOS=8\relax}%
\@namedef{dgo@9}{\dg@LBLPOS=9\relax}%
\def\dgt@e{\dg@DX=\@ne \dg@DY=\z@ \dg@SIZE=\@ne}%
\def\dgt@w{\dg@DX=\m@ne \dg@DY=\z@ \dg@SIZE=\@ne}%
\def\dgt@n{\dg@DX=\z@ \dg@DY=\@ne \dg@SIZE=\@ne}%
\def\dgt@s{\dg@DX=\z@ \dg@DY=\m@ne \dg@SIZE=\@ne}%
\def\dgt@ne{\dg@DX=\@ne \dg@DY=\@ne \dg@SIZE=\@ne}%
\def\dgt@se{\dg@DX=\@ne \dg@DY=\m@ne \dg@SIZE=\@ne}%
\def\dgt@nw{\dg@DX=\m@ne \dg@DY=\@ne \dg@SIZE=\@ne}%
\def\dgt@sw{\dg@DX=\m@ne \dg@DY=\m@ne \dg@SIZE=\@ne}%
\def\dgt@nne{\dg@DX=\@ne \dg@DY=\tw@ \dg@SIZE=\@ne}%
\def\dgt@nnw{\dg@DX=\m@ne \dg@DY=\tw@ \dg@SIZE=\@ne}%
\def\dgt@sse{\dg@DX=\@ne \dg@DY=-\tw@ \dg@SIZE=\@ne}%
\def\dgt@ssw{\dg@DX=\m@ne \dg@DY=-\tw@ \dg@SIZE=\@ne}%
\def\dgt@ene{\dg@DX=\tw@ \dg@DY=\@ne \dg@SIZE=\tw@}%
\def\dgt@ese{\dg@DX=\tw@ \dg@DY=\m@ne \dg@SIZE=\tw@}%
\def\dgt@wnw{\dg@DX=-\tw@ \dg@DY=\@ne \dg@SIZE=\tw@}%
\def\dgt@wsw{\dg@DX=-\tw@ \dg@DY=\m@ne \dg@SIZE=\tw@}%
\def\dggeometry{
   \dg@ZTEMP=\dg@XGRID \multiply\dg@ZTEMP\tw@
   \ifnum\dg@YGRID=\z@ \dg@ZTEMP=\tw@
   \else \divide\dg@ZTEMP\dg@YGRID \fi
   \ifnum\dg@ZTEMP>\f@ur \dg@ZTEMP=\f@ur \fi
   \ifnum\dg@ZTEMP<\@ne \dg@ZTEMP=\@ne \fi
   \unitlength=2sp\relax
   \ifnum\dg@ZTEMP<\tw@
      \advance\dg@ZTEMP\@ne
      \multiply\unitlength\dg@YGRID
   \else
      \multiply\unitlength\dg@XGRID \divide\unitlength\dg@ZTEMP
   \fi
   \dg@XGRID=\dg@ZTEMP \dg@YGRID=\tw@
   \dg@rmcommondiv\tw@\dg@XGRID\dg@YGRID
   \divide\unitlength\dg@YGRID \divide\unitlength\@m\relax}
\@namedef{dgo@..}{\let\dg@VECTOR=\dg@dotvector}%

\def\dg@dotvector(#1,#2)#3{%
   \begingroup
   \dg@XTEMP=#1\relax \dg@YTEMP=#2\relax
   \let\dg@NDOTS=\dg@XEND \let\dg@DOTDIAM=\dg@WEND
   \dg@NDOTS=\unitlength \multiply\dg@NDOTS #3\relax
   \dg@ZTEMP=\dg@YTEMP \dg@changesign\dg@YTEMP\dg@ZTEMP
   \ifnum\dg@XTEMP>\z@
      \ifnum\dg@YTEMP>\dg@XTEMP
         \multiply\dg@NDOTS\dg@YTEMP \divide\dg@NDOTS\dg@XTEMP \fi
   \else\ifnum\dg@XTEMP<\z@
      \ifnum\dg@YTEMP>-\dg@XTEMP
         \multiply\dg@NDOTS\dg@YTEMP \divide\dg@NDOTS-\dg@XTEMP \fi
   \fi\fi
   \dg@YTEMP=\dg@ZTEMP
   \divide\dg@NDOTS\dgDOTSPACING
   \ifnum\dg@NDOTS>\z@\else \dg@NDOTS=\@ne \fi
   \dg@ZTEMP=\unitlength \multiply\dg@ZTEMP #3\relax
   \divide\dg@ZTEMP\dg@NDOTS
   \ifnum\dg@XTEMP=\z@
      \dg@changesign\dg@ZTEMP\dg@YTEMP \dg@YTEMP=\dg@ZTEMP
   \else
      \dg@changesign\dg@ZTEMP\dg@XTEMP
      \multiply\dg@YTEMP\dg@ZTEMP \divide\dg@YTEMP\dg@XTEMP
      \dg@XTEMP=\dg@ZTEMP
   \fi
   \divide\dg@XTEMP\unitlength \divide\dg@YTEMP\unitlength
   \begin{picture}(0,0)
      \dg@DOTDIAM=\dgDOTSIZE \divide\dg@DOTDIAM\unitlength
      \multiput(0,0)(\dg@XTEMP,\dg@YTEMP){\dg@NDOTS}{%
         \circle*{\dg@DOTDIAM}}%
      \multiply\dg@XTEMP\dg@NDOTS \multiply\dg@YTEMP\dg@NDOTS
      \put(\dg@XTEMP,\dg@YTEMP){\vector(#1,#2){0}}%
   \end{picture}%
   \endgroup}%
\ifx\@auxdefsloaded\relax
\else\let\@auxdefsloaded\relax

\def\newenvironment{%
   \@ifnextchar *{\@@newenv{\global\@ignoretrue}}{\@@newenv{}*}}

\def\@@newenv#1*#2{%
   \@ifnextchar [{\@newenv{#1}{#2}}{\@newenv{#1}{#2}[0]}}

\long\def\@newenv#1#2[#3]#4#5{%
   \expandafter\newcommand\csname#2\endcsname[#3]{#4}%
   \expandafter\long\expandafter\def\csname end#2\endcsname{#5#1}}

\def\renewenvironment{%
   \@ifnextchar *{%
      \@@renewenv{\global\@ignoretrue}}{\@@renewenv{}*}}

\def\@@renewenv#1*#2{%
   \@ifnextchar [{\@renewenv{#1}{#2}}{\@renewenv{#1}{#2}[0]}}

\long\def\@renewenv#1#2[#3]#4#5{%
   \expandafter\renewcommand\csname#2\endcsname[#3]{#4}%
   \expandafter\long\expandafter\def\csname end#2\endcsname{#5#1}}

\def\newoptcommand#1#2{%
   \@ifnextchar [{\@optargdef#1#2}{\@optargdef#1#2[1]}}

\def\renewoptcommand#1#2{%
   \edef\@tempa{\expandafter\@cdr\string#1\@nil}%
   \@ifundefined{\@tempa}{%
      \@latexerr{\string#1\space undefined}\@ehc}{}%
   \@ifnextchar [{\@reoptargdef#1#2}{\@reoptargdef#1#2[1]}}

\long\def\@optargdef#1#2[#3]#4{%
   \@ifdefinable #1{\@reoptargdef#1#2[#3]{#4}}}

\long\def\@reoptargdef#1#2[#3]#4{%
   \@tempcnta#3\relax \@tempcntb \@ne
   \let#1\relax \let\@tempa\relax
   \edef\@tempb{\long\def\csname\string#1\endcsname
      [\@tempa\the\@tempcntb]}%
   \advance\@tempcntb \@ne \advance\@tempcnta \m@ne
   \@whilenum\@tempcnta>0\do{%
      \edef\@tempb{\@tempb\@tempa\the\@tempcntb}%
      \advance\@tempcntb \@ne \advance\@tempcnta \m@ne}%
   \let\@tempa=##\@tempb{#4}%
   \def#1{\@ifnextchar [{\csname\string#1\endcsname}{%
      \csname\string#1\endcsname[#2]}}}

\def\newoptenvironment{%
   \@ifnextchar *{\@@newoptenv{\global\@ignoretrue}}{%
      \@@newoptenv{}*}}

\def\@@newoptenv#1*#2#3{%
   \@ifnextchar [{\@newoptenv{#1}{#2}{#3}}{%
      \@newoptenv{#1}{#2}{#3}[0]}}

\long\def\@newoptenv#1#2#3[#4]#5#6{%
   \expandafter\newoptcommand\csname#2\endcsname{#3}[#4]{#5}%
   \expandafter\long\expandafter\def\csname end#2\endcsname{#6#1}}

\def\renewoptenvironment{%
   \@ifnextchar *{\@@renewoptenv{\global\@ignoretrue}}{%
      \@@renewoptenv{}*}}

\def\@@renewoptenv#1*#2#3{%
   \@ifnextchar [{\@renewoptenv{#1}{#2}{#3}}{%
      \@renewoptenv{#1}{#2}{#3}[0]}}

\long\def\@renewoptenv#1#2#3[#4]#5#6{%
   \expandafter\renewoptcommand\csname#2\endcsname{#3}[#4]{#5}%
   \expandafter\long\expandafter\def\csname end#2\endcsname{#6#1}}

\newcounter{keepoptional}
\newcounter{optnestctr}

\newoptenvironment*{optional}{1}[1]{%
   \ifnum#1>\value{keepoptional}\relax
      \setcounter{optnestctr}{0}\@powerup\expandafter\@endopt\fi
   \ignorespaces}{}

\def\@powerup{\catcode`\{=12 \catcode`\}=12 \catcode`\\=12 \relax}
\def\@powerdown{\catcode`\{=1 \catcode`\}=2 \catcode`\\=0 \relax}
\begingroup \catcode`|=0 \catcode`[=1 \catcode`]=2 \@powerup
   |long|gdef|@endopt#1\end{optional}[%
      |@beginopt#1\begin{optional}*]%
   |long|gdef|@beginopt#1\begin{optional}[%
      |@ifstar[%
         |ifnum|value[optnestctr]>0|relax
            |addtocounter[optnestctr][-1]|let|@zzzap=|@endopt
         |else
            |@powerdown|end[optional]|let|@zzzap=|relax|fi
	 |@zzzap]%
      [%else
         |addtocounter[optnestctr][1]|@beginopt]]%
|endgroup

\fi

\newcount\dg@HORIZ      \newcount\dg@VERT
\newcount\dg@XLPAD      \newcount\dg@YBPAD
\newcount\dg@XRPAD      \newcount\dg@YTPAD
\newcount\dg@X          \newcount\dg@Y
\newcount\dg@XGRID      \newcount\dg@YGRID
\newcount\dg@SIZE
\newcount\dg@USERSIZE
\newcount\dg@DX         \newcount\dg@DY
\newcount\dg@XLBL       \newcount\dg@YLBL
\newcount\dg@XOFFSET    \newcount\dg@YOFFSET
\newcount\dg@LBLOFF
\newcount\dg@XLBLOFF    \newcount\dg@YLBLOFF
\newcount\dg@LBLPOS
\newcount\dg@XTEMP      \newcount\dg@YTEMP
\newcount\dg@ZTEMP
\newcount\dg@XNODE      \newcount\dg@YNODE
\newcount\dg@XEND       \newcount\dg@YEND
\newcount\dg@WEND       \newcount\dg@HEND
\newcount\dg@COUNT
\newbox\dg@NODEBOX
\newoptenvironment*{diagram}{}[1]{%
   \global\advance\dg@COUNT\@ne \typeout{[diagram \the\dg@COUNT:}%
   \let\node=\dg@node \let\\=\dg@cr \let\arrow=\dg@arrow
   \def\dg@BIGNODE{#1}%
   \ifx\dg@BIGNODE\@empty\else
      \setbox\dg@NODEBOX=\hbox{$\dgeverynode{#1}$}%
      \dg@XTEMP=\wd\dg@NODEBOX
      \dg@YTEMP=\ht\dg@NODEBOX \advance\dg@YTEMP\dp\dg@NODEBOX
      \ifnum\dg@YTEMP=\z@ \dg@ZTEMP=\z@
      \else \dg@ZTEMP=\dg@XTEMP \divide\dg@ZTEMP\dg@YTEMP\fi
      \ifnum\dg@ZTEMP>\dgMAXSQUARE
         \ifnum\dg@ZTEMP<\dgMINDBLSQ
		 	\dg@XGRID=\thr@@ \dg@YGRID=\tw@
         \else \dg@XGRID=\tw@ \dg@YGRID=\@ne \fi
      \else \dg@XGRID=\@ne \dg@YGRID=\@ne \fi
      \advance\dg@XTEMP\dgHORIZPAD \advance\dg@YTEMP\dgVERTPAD
      \dg@XLPAD=\dg@XTEMP \divide\dg@XLPAD\tw@
      \dg@XRPAD=\dg@XTEMP \divide\dg@XRPAD\tw@
      \dg@YBPAD=\dg@YTEMP \divide\dg@YBPAD\tw@
      \dg@YTPAD=\dg@YTEMP \divide\dg@YTPAD\tw@
      \advance\dg@XTEMP\dgARROWLENGTH
      \divide\dg@XTEMP\dg@XGRID \divide\dg@XTEMP\@m
      \unitlength=1sp\relax \multiply\unitlength\dg@XTEMP
   \fi
   \typeout{saving...}%
   \dg@X=-\@ne \dg@Y=\z@ \dg@HORIZ=\z@\relax
   \let\dg@SLIST=\@empty
   \let\dg@NLIST=\@empty \let\dg@ALIST=\@empty
   \let\dg@PASS=\dg@savepass
}{%\enddiagram
   \dg@VERT=-\dg@Y\relax
   \typeout{calculating...}%
   \ifx\dg@BIGNODE\@empty
      \dg@XGRID=\z@ \dg@YGRID=\z@
      \dg@XLPAD=\z@ \dg@XRPAD=\z@ \dg@YBPAD=\z@ \dg@YTPAD=\z@
      \let\dg@PASS=\dg@geompass
      \dg@NLIST \dg@ALIST
      \ifnum\dg@XGRID>\z@\else \dg@XGRID=\dgHORIZPAD \fi
      \ifnum\dg@YGRID>\z@\else \dg@YGRID=\dgVERTPAD \fi
      \dggeometry
      \ifdim\unitlength>\z@\else
         \unitlength=\dgARROWLENGTH \divide\unitlength\@m \fi
      \ifnum\dg@XGRID>\z@\else \dg@XGRID=\@ne \fi
      \ifnum\dg@YGRID>\z@\else \dg@YGRID=\@ne \fi
   \fi
   \dg@LBLOFF=\dgLABELOFFSET \divide\dg@LBLOFF\unitlength
   \multiply\dg@HORIZ\@m \multiply\dg@HORIZ\dg@XGRID
   \multiply\dg@VERT\@m \multiply\dg@VERT\dg@YGRID
   \divide\dg@XLPAD\unitlength \divide\dg@XRPAD\unitlength
   \divide\dg@YBPAD\unitlength \divide\dg@YTPAD\unitlength
   \typeout{drawing...}%
   \let\dg@PASS=\dg@drawpass
   \hspace*{\dg@XLPAD\unitlength}%
   \vcenter{%
      \hsize=0pt\relax\parindent=0pt\relax
      \parskip=0pt\relax\baselineskip=0pt\relax
      \vspace*{\dg@YTPAD\unitlength}%
      \begin{picture}(0,0)(0,0)\dg@NLIST\dg@ALIST\end{picture}%
      \raisebox{\z@}[\z@][\dg@VERT\unitlength]{}%
      \vspace*{\dg@YBPAD\unitlength}%
   }
   \hspace*{\dg@HORIZ\unitlength}%
   \hspace*{\dg@XRPAD\unitlength}%
   \typeout{done]}}%

\def\dg@savepass{s}
\def\dg@geompass{g}
\def\dg@drawpass{d}

\newoptcommand{\dg@node}{\@ne}[2]{%
   \ifx\dg@PASS\dg@savepass
      \dg@XTEMP=#1\relax
      \ifnum\dg@XTEMP<\@ne \dg@XTEMP=\@ne\fi
      \advance\dg@X\dg@XTEMP
      \ifnum\dg@HORIZ<\dg@X \dg@HORIZ=\dg@X \fi
      \setbox\dg@NODEBOX=\hbox{$\dgeverynode{#2}$}%
      \dg@XTEMP=\wd\dg@NODEBOX \advance\dg@XTEMP\dgHORIZPAD
      \dg@YTEMP=\ht\dg@NODEBOX \advance\dg@YTEMP\dp\dg@NODEBOX
      \advance\dg@YTEMP\dgVERTPAD
      \toks\z@=\expandafter{\dg@SLIST}%
      \edef\dg@SLIST{\the\toks\z@
         ,\noexpand\dg@XNODE=\number\dg@X\noexpand\relax
         \noexpand\dg@YNODE=\number\dg@Y\noexpand\relax
         \noexpand\dg@XTEMP=\number\dg@XTEMP\noexpand\relax
         \noexpand\dg@YTEMP=\number\dg@YTEMP\noexpand\relax}%
      \toks\z@=\expandafter{\dg@NLIST}%
      \toks\tw@={\dg@node{#2}}%
      \edef\dg@NLIST{\the\toks\z@
         \noexpand\dg@X=\number\dg@X\noexpand\relax
         \noexpand\dg@Y=\number\dg@Y\noexpand\relax
         \the\toks\tw@}%
   \else\ifx\dg@PASS\dg@geompass
      \ifnum\dg@X=\z@
         \dg@getnodesize
            {\dg@SLIST}{\dg@X}{\dg@Y}{\dg@WEND}{\dg@HEND}%
         \divide\dg@WEND\tw@
         \ifnum\dg@XLPAD<\dg@WEND \dg@XLPAD=\dg@WEND \fi\fi
      \ifnum\dg@X=\dg@HORIZ
         \dg@getnodesize
            {\dg@SLIST}{\dg@X}{\dg@Y}{\dg@WEND}{\dg@HEND}%
         \divide\dg@WEND\tw@
         \ifnum\dg@XRPAD<\dg@WEND \dg@XRPAD=\dg@WEND \fi\fi
      \ifnum\dg@Y=\z@
         \dg@getnodesize
            {\dg@SLIST}{\dg@X}{\dg@Y}{\dg@WEND}{\dg@HEND}%
         \divide\dg@HEND\tw@
         \ifnum\dg@YTPAD<\dg@HEND \dg@YTPAD=\dg@HEND \fi\fi
      \ifnum\dg@Y=-\dg@VERT\relax
         \dg@getnodesize
            {\dg@SLIST}{\dg@X}{\dg@Y}{\dg@WEND}{\dg@HEND}%
         \divide\dg@HEND\tw@
         \ifnum\dg@YBPAD<\dg@HEND \dg@YBPAD=\dg@HEND \fi\fi
   \else\ifx\dg@PASS\dg@drawpass
      \dg@XNODE=\dg@X \multiply\dg@XNODE\@m
      \multiply\dg@XNODE\dg@XGRID
      \dg@YNODE=\dg@Y \multiply\dg@YNODE\@m
      \multiply\dg@YNODE\dg@YGRID
      \setbox\dg@NODEBOX=\hbox{$\dgeverynode{#2}$}%
      \put(\dg@XNODE,\dg@YNODE){\dg@makebox{\box\dg@NODEBOX}}%
   \fi\fi\fi}%

\newoptcommand{\dg@cr}{\@ne}[1]{%
   \ifx\dg@PASS\dg@savepass
      \dg@YTEMP=#1\relax
      \ifnum\dg@YTEMP<\@ne \dg@YTEMP=\@ne \fi
      \advance\dg@Y -\dg@YTEMP\relax
      \dg@X=-\@ne\relax\fi}%
\newoptcommand{\dg@arrow}{\@ne}[2]{%
   \begingroup
   \dg@USERSIZE=#1\relax
   \ifnum\dg@USERSIZE<\@ne \dg@USERSIZE=\@ne \fi
   \dg@parse{#2}%
   \ifx\dg@PASS\dg@savepass
      \ifx\dg@label\dgl@b \let\dg@label=\dgl@t \fi
      \ifx\dg@label\dgl@r \let\dg@label=\dgl@l \fi
      \let\dg@process=\dg@save
   \else\ifx\dg@PASS\dg@geompass
      \let\dg@process=\dg@ignore
      \dg@geomcalc
   \else\ifx\dg@PASS\dg@drawpass
      \let\dg@process=\dg@draw
      \dg@drawcalc
   \fi\fi\fi
   \dg@label{\dg@process{#1}{#2}}}%

\newoptcommand{\arrow}{\@ne}[2]{%
   \dg@parse{#2}%
   \ifx\dg@label\dgl@b \let\dg@label=\dgl@t \fi
   \ifx\dg@label\dgl@r \let\dg@label=\dgl@l \fi
   \dg@label{\dg@textarrow{#1}{#2}}}%

\def\dg@textarrow#1#2#3#4{%
   \mathop{{\dgHORIZPAD=0pt\relax\dgVERTPAD=0pt\relax
      \begin{diagram}
         \node{}\arrow[#1]{#2}{#3}{#4}\node{}
      \end{diagram}}}}

\def\dg@parse#1{%
   \let\dg@label=\dgl@ \dgo@
   \let\dg@type=\@empty \let\dg@lbltype=\@empty
   \@for\dg@list:=#1\do{%
      \ifx\dg@type\@empty \let\dg@type=\dg@list
      \else\ifx\dg@lbltype\@empty \let\dg@lbltype=\dg@list
         \@ifundefined{dgo@\dg@list}{}{\@nameuse{dgo@\dg@list}}%
      \else
         \@ifundefined{dgo@\dg@list}{}{\@nameuse{dgo@\dg@list}}%
      \fi\fi}%
   \@ifundefined{dgt@\dg@type}{\dgt@e}{\@nameuse{dgt@\dg@type}}%
   \@ifundefined{dgl@\dg@lbltype}{}{%
      \dg@letname\dg@label{dgl@\dg@lbltype}}}

\def\dg@draw#1#2#3#4{%
   \put(\dg@X,\dg@Y){\dg@makebox{%
      \begin{picture}(0,0)%
         \thinlines
         \put(\dg@XOFFSET,\dg@YOFFSET){%
            \dg@VECTOR(\dg@DX,\dg@DY){\dg@SIZE}}%
         \put(\dg@XLBL,\dg@YLBL){\dg@makebox{%
            \begin{picture}(0,0)%
               \put(\dg@XLBLOFF,\dg@YLBLOFF){%
                  \dg@makebox[\dg@LBLONE]{$\dgeverylabel{#3}$}}%
               \put(-\dg@XLBLOFF,-\dg@YLBLOFF){%
                  \dg@makebox[\dg@LBLTWO]{$\dgeverylabel{#4}$}}%
            \end{picture}}}%
      \end{picture}}}%
   \endgroup}%

\def\dg@save#1#2#3#4{%
   \endgroup % to match \dg@arrow's \begingroup
   \toks\z@=\expandafter{\dg@ALIST}%
   \toks\tw@={\dg@arrow[#1]{#2}{#3}{#4}}%
   \edef\dg@ALIST{\the\toks\z@%
      \noexpand\dg@X=\number\dg@X\noexpand\relax
      \noexpand\dg@Y=\number\dg@Y\noexpand\relax
      \the\toks\tw@}}%

\def\dg@ignore#1#2#3#4{\endgroup}

\def\dg@geomcalc{%
   \dg@XEND=\dg@SIZE \multiply\dg@XEND\dg@USERSIZE
   \ifnum\dg@DX=\z@
      \dg@YEND=\dg@XEND \dg@XEND=\z@
      \dg@changesign\dg@YEND\dg@DY
   \else
      \dg@changesign\dg@XEND\dg@DX \dg@YEND=\dg@XEND
      \multiply\dg@YEND\dg@DY \divide\dg@YEND\dg@DX
   \fi
   \advance\dg@XEND\dg@X \advance\dg@YEND\dg@Y
   \dg@getnodesize
      {\dg@SLIST}{\dg@XEND}{\dg@YEND}{\dg@WEND}{\dg@HEND}%
   \dg@XOFFSET=\dg@WEND \dg@YOFFSET=\dg@HEND
   \dg@getnodesize
      {\dg@SLIST}{\dg@X}{\dg@Y}{\dg@WEND}{\dg@HEND}%
   \advance\dg@XOFFSET\dg@WEND \divide\dg@XOFFSET\tw@
   \advance\dg@YOFFSET\dg@HEND \divide\dg@YOFFSET\tw@
   \dg@XTEMP=\dgARROWLENGTH \dg@YTEMP=\dgARROWLENGTH
   \ifnum\dg@DX>\z@
      \dg@ZTEMP=\dg@DX \multiply\dg@XTEMP\dg@DX
   \else \dg@ZTEMP=-\dg@DX \multiply\dg@XTEMP -\dg@DX \fi
   \ifnum\dg@DY>\z@
      \advance\dg@ZTEMP\dg@DY \multiply\dg@YTEMP\dg@DY
   \else \advance\dg@ZTEMP -\dg@DY \multiply\dg@YTEMP -\dg@DY\fi
   \ifnum\dg@ZTEMP=\z@\else
      \divide\dg@XTEMP\dg@ZTEMP \divide\dg@YTEMP\dg@ZTEMP
      \advance\dg@XOFFSET\dg@XTEMP \advance\dg@YOFFSET\dg@YTEMP
   \fi
   \divide\dg@XOFFSET\dg@SIZE \divide\dg@XOFFSET\dg@USERSIZE
   \divide\dg@YOFFSET\dg@SIZE \divide\dg@YOFFSET\dg@USERSIZE
   \ifnum\dg@DX=\z@ \dg@XOFFSET=\z@ \fi
   \ifnum\dg@DY=\z@ \dg@YOFFSET=\z@ \fi
   \ifnum\dg@XGRID<\dg@XOFFSET \global\dg@XGRID=\dg@XOFFSET\fi
   \ifnum\dg@YGRID<\dg@YOFFSET \global\dg@YGRID=\dg@YOFFSET\fi
   \relax}

\def\dg@drawcalc{%
   \dg@XEND=\dg@SIZE \multiply\dg@XEND\dg@USERSIZE
   \ifnum\dg@DX=\z@
      \dg@YEND=\dg@XEND \dg@XEND=\z@
      \dg@changesign\dg@YEND\dg@DY
   \else
      \dg@changesign\dg@XEND\dg@DX \dg@YEND=\dg@XEND
      \multiply\dg@YEND\dg@DY \divide\dg@YEND\dg@DX
   \fi
   \advance\dg@XEND\dg@X \advance\dg@YEND\dg@Y
   \dg@getnodesize
      {\dg@SLIST}{\dg@XEND}{\dg@YEND}{\dg@WEND}{\dg@HEND}%
   \divide\dg@WEND\unitlength \divide\dg@HEND\unitlength
   \multiply\dg@DX\dg@XGRID \multiply\dg@DY\dg@YGRID
   \dg@rmcommondiv\tw@\dg@DX\dg@DY
   \dg@rmcommondiv\tw@\dg@DX\dg@DY %[sic]
   \dg@rmcommondiv\thr@@\dg@DX\dg@DY
   \multiply\dg@SIZE\dg@USERSIZE \multiply\dg@SIZE\@m
   \ifnum\dg@DX=\z@
      \multiply\dg@SIZE\dg@YGRID
      \divide\dg@HEND\tw@ \advance\dg@SIZE -\dg@HEND
      \dg@getnodesize
         {\dg@SLIST}{\dg@X}{\dg@Y}{\dg@WEND}{\dg@YOFFSET}%
      \divide\dg@YOFFSET\unitlength \divide\dg@YOFFSET\tw@
      \advance\dg@SIZE -\dg@YOFFSET
      \dg@XOFFSET=\z@
      \def\dg@LBLONE{r}\def\dg@LBLTWO{l}%
      \dg@XLBL=\z@ \dg@YLBL=\dg@SIZE
      \multiply\dg@YLBL\dg@LBLPOS
      \divide\dg@YLBL\dgARROWPARTS\relax
      \advance\dg@YLBL\dg@YOFFSET
      \dg@changesign\dg@YLBL\dg@DY
      \dg@changesign\dg@YOFFSET\dg@DY
   \else
      \multiply\dg@SIZE\dg@XGRID
      \ifnum\dg@DY=\z@
         \divide\dg@WEND\tw@ \advance\dg@SIZE -\dg@WEND
         \dg@getnodesize
            {\dg@SLIST}{\dg@X}{\dg@Y}{\dg@XOFFSET}{\dg@HEND}%
         \divide\dg@XOFFSET\unitlength \divide\dg@XOFFSET\tw@
         \advance\dg@SIZE -\dg@XOFFSET
         \dg@YOFFSET=\z@
         \def\dg@LBLONE{b}\def\dg@LBLTWO{t}%
         \dg@YLBL=\z@ \dg@XLBL=\dg@SIZE
         \multiply\dg@XLBL\dg@LBLPOS
         \divide\dg@XLBL\dgARROWPARTS\relax
         \advance\dg@XLBL\dg@XOFFSET
         \dg@changesign\dg@XLBL\dg@DX
         \dg@changesign\dg@XOFFSET\dg@DX
      \else
         \divide\dg@WEND\tw@ \divide\dg@HEND\tw@
         \multiply\dg@HEND\dg@DX \divide\dg@HEND\dg@DY
         \ifnum\dg@HEND<\z@ \multiply\dg@HEND\m@ne \fi
         \ifnum\dg@WEND<\dg@HEND \advance\dg@SIZE -\dg@WEND
         \else \advance\dg@SIZE -\dg@HEND \fi
         \dg@getnodesize
            {\dg@SLIST}{\dg@X}{\dg@Y}{\dg@WEND}{\dg@HEND}%
         \divide\dg@WEND\unitlength \divide\dg@WEND\tw@
         \divide\dg@HEND\unitlength \divide\dg@HEND\tw@
         \multiply\dg@HEND\dg@DX \divide\dg@HEND\dg@DY
         \ifnum\dg@HEND<\z@ \multiply\dg@HEND\m@ne \fi
         \ifnum\dg@WEND<\dg@HEND \dg@XOFFSET=\dg@WEND
         \else \dg@XOFFSET=\dg@HEND \fi
         \advance\dg@SIZE -\dg@XOFFSET
         \dg@changesign\dg@XOFFSET\dg@DX
         \dg@YOFFSET=\dg@XOFFSET
         \multiply\dg@YOFFSET\dg@DY \divide\dg@YOFFSET\dg@DX
         \def\dg@LBLONE{br}\def\dg@LBLTWO{tl}%
         \ifnum\dg@DX<\z@ \ifnum\dg@DY>\z@
            \def\dg@LBLONE{bl}\def\dg@LBLTWO{tr}\fi\fi
         \ifnum\dg@DX>\z@ \ifnum\dg@DY<\z@
            \def\dg@LBLONE{bl}\def\dg@LBLTWO{tr}\fi\fi
         \dg@XLBL=\dg@SIZE
         \multiply\dg@XLBL\dg@LBLPOS
         \divide\dg@XLBL\dgARROWPARTS\relax
         \dg@changesign\dg@XLBL\dg@DX
         \dg@YLBL=\dg@XLBL
         \multiply\dg@YLBL\dg@DY \divide\dg@YLBL\dg@DX
         \advance\dg@XLBL\dg@XOFFSET
         \advance\dg@YLBL\dg@YOFFSET
      \fi
   \fi
   \dg@XLBLOFF=-\dg@DY \dg@changesign\dg@XLBLOFF\dg@DX
   \dg@YLBLOFF=\dg@DX \dg@changesign\dg@YLBLOFF\dg@DX
   \ifnum\dg@DX=\z@ \dg@XLBLOFF=\m@ne \fi
   \dg@XTEMP=\dg@DX \dg@changesign\dg@XTEMP\dg@DX
   \dg@YTEMP=\dg@DY \dg@changesign\dg@YTEMP\dg@DY
   \ifnum\dg@YTEMP>\dg@XTEMP \dg@XTEMP=\dg@YTEMP \fi
   \ifnum\dg@XTEMP=\z@ \dg@XTEMP=\@ne \fi
   \multiply\dg@XLBLOFF\dg@LBLOFF \divide\dg@XLBLOFF\dg@XTEMP
   \multiply\dg@YLBLOFF\dg@LBLOFF \divide\dg@YLBLOFF\dg@XTEMP
   \multiply\dg@X\@m \multiply\dg@X\dg@XGRID
   \multiply\dg@Y\@m \multiply\dg@Y\dg@YGRID
   \relax}%

\def\dg@rmcommondiv#1#2#3{%
   \dg@XTEMP=#2\relax
   \divide\dg@XTEMP #1\relax \multiply\dg@XTEMP #1\relax
   \dg@YTEMP=#3\relax
   \divide\dg@YTEMP #1\relax \multiply\dg@YTEMP #1\relax
   \ifnum\dg@XTEMP=#2\relax \ifnum\dg@YTEMP=#3\relax
      \divide#2#1\relax \divide#3#1\relax \fi\fi}%

\def\dg@changesign#1#2{%
   \ifnum #2<\z@ \multiply#1\m@ne
   \else\ifnum #2=\z@ #1=\z@ \fi\fi}%

\def\dg@getnodesize#1#2#3#4#5{%
   #4=\z@\relax #5=\z@\relax
   \expandafter\@for\expandafter\dg@trynode
   \expandafter:\expandafter=#1\do{%
      \dg@XNODE=\m@ne % impossible (in case \dg@trynode is empty)
      \dg@trynode
      \ifnum #2=\dg@XNODE \ifnum #3=\dg@YNODE
         #4=\dg@XTEMP\relax #5=\dg@YTEMP\relax\fi\fi}}%

\long\def\dg@sinkbaseline#1{%
   \leavevmode\hbox{\vbox{%
      \lineskiplimit=16383pt\relax\lineskip=0pt\relax
      \baselineskip=-1000pt\relax
      \parindent=0pt\relax\parskip=0pt\relax
      \hbox{#1}\rule{0pt}{0pt}}}}%
\newoptcommand{\dg@makebox}{}[2]{\dg@sinkbaseline{%
   \expandafter\makebox\expandafter(\expandafter
      0\expandafter,\expandafter0\expandafter)\expandafter
      [#1]{#2}}}%

\def\dg@novector(#1,#2)#3{}%

\def\dg@letname#1#2{%
   \relax\expandafter
   \let\expandafter #1\csname #2\endcsname\relax}%

\def\dgl@#1{#1{}{}}%
\def\dgl@t#1#2{#1{#2}{}}%
\def\dgl@b#1#2{#1{}{#2}}%
\def\dgl@tb#1#2#3{#1{#2}{#3}}%
\def\dgl@l#1#2{#1{#2}{}}%
\def\dgl@r#1#2{#1{}{#2}}%
\def\dgl@lr#1#2#3{#1{#2}{#3}}%
\makeatother

\expandafter\ifx\csname >amspapermacs\endcsname\relax
\expandafter\def\csname >amspapermacs\endcsname{done}
\else\endinput\fi
\newtheorem{thm}{Theorem}[section]
\newtheorem{cor}[thm]{Corollary}
\newtheorem{lem}[thm]{Lemma}
\newtheorem{prop}[thm]{Proposition}

\theoremstyle{definition}
\newtheorem{definition}[thm]{Definition}

\theoremstyle{remark}
\newtheorem{remark}[thm]{Remark}

\numberwithin{equation}{section}

\newcommand{\thmref}[1]{Theorem~\ref{#1}}
\newcommand{\secref}[1]{Section~\ref{#1}}
\newcommand{\lemref}[1]{Lemma~\ref{#1}}
\newcommand{\propref}[1]{Prop\-o\-si\-tion~\ref{#1}}

\makeatletter
\newif\ifdraft
\draftfalse

\let\foob@r=\label
\let\f@@index\index
\def\t@index#1{\relax}
\def\label#1{\foob@r{#1}\ifdraft
\index{\string\labelentry{\@currenvir}{#1}}
\ifinner
\else\leavevmode\marginpar[\hfill
{\tiny\normalshape #1}]{{\tiny\normalshape
#1}}\fi\fi\ignorespaces}
\def\@lem@#1{Lemma~\ref{#1}}
\def\@equation@#1{Equation~(\ref{#1})}
\let\@align@=\@equation@
\let\@gather@=\@equation@
\def\@thm@#1{Theorem~\ref{#1}}
\def\@cor@#1{Corollary~\ref{#1}}
\def\@prop@#1{Proposition~\ref{#1}}
\def\@definition@#1{Definition~\ref{#1}}
\def\@remark@#1{Remark~\ref{#1}}
\def\@ex@#1{Example~\ref{#1}}
\def\dotfill{\leaders\hbox to0.5em{\hss.\hss}\hfill}
\newdimen\lcrskip
\def\labelentry#1#2#3{\vskip0pt plus 0.25pt\hbox
to\hsize{\hskip\lcrskip\@ifundefined{@#1@}{{\tt{\uppercase{#1}}}~\ref{#2}}
{\csname @#1@\endcsname{#2}}\dotfill{#2}\quad
(Page #3)\hskip\lcrskip}}
\def\labeltable{\ifdraft\clearpage\immediate\closeout\@indexfile
\section*{Label Cross References}
\input \jobname.idx
\else\relax\fi}
\def\gmeasure@#1{\let\f@@index\index\let\index\t@index\gwidth@
\z@\gmaxwidth@\z@\setbox@ne\vbox{\Let@
 \firstchoice@false\let\tag\gobble@tag
 \halign{\setboxz@h{$\m@th\displaystyle{##}$}\global\gwidth@\wdz@
\ifdim\gwidth@>\gmaxwidth@\global\gmaxwidth@\gwidth@\fi
 &\@gobble{##}\crcr#1\crcr}}\let\index\f@@index}
\def\measure@#1{\let\f@@index\index\let\index\t@index
\lwidth@\z@\rwidth@\z@\maxlwidth@\z@\maxrwidth@\z@
 \global\and@\z@
 \setbox@ne\vbox{%
   \everycr{\noalign{\global\tag@false\global\and@\z@}}\Let@
   \let\tag\gobble@tag
   \let\notag\@empty \let\nonumber\@empty
   \firstchoice@false
    \halign{\setboxz@h{$\m@th\displaystyle{\@lign##}$}%
     \global\lwidth@\wdz@
     \ifdim\lwidth@>\maxlwidth@\global\maxlwidth@\lwidth@\fi
     \global\advance\and@\@ne
     &\setboxz@h{$\m@th\displaystyle{{}\@lign##}$}%
     \global\rwidth@\wdz@
     \ifdim\rwidth@>\maxrwidth@\global\maxrwidth@\rwidth@\fi
     \global\advance\and@\@ne
     &\Tag@\@gobble{##}\crcr#1\crcr}}%
 \totwidth@\maxlwidth@\advance\totwidth@\maxrwidth@
\let\index\f@@index}
\newbox\di@b@x

\newbox\di@b@x
\newenvironment{ediagram*}{\setbox\di@b@x=\hbox\bgroup$
\begin{diagram}}{\end{diagram}$\egroup%
\begin{equation*}\copy\di@b@x\end{equation*}\global\@ignoretrue}
\makeatother

\newcount\hours
\newcount\minutes       %   For computing the time of day on
\def\timeofday{%    Must be computed when called if preloaded
\hours=\time
\minutes=\hours
\divide\hours by60
\multiply\hours by60
\advance\minutes by-\hours
\divide\hours by60
\ifnum\hours>9\else0\fi\the\hours:\ifnum\minutes>9\else
0\fi\the\minutes}
\def\predate{\date{\the\day\ \ifcase\month\or
  January\or February\or March\or April\or May\or June\or July\or
	August\or September\or October\or November\or
	   December\fi\ \the\year\ --- \timeofday\ --- Preliminary
		  Version}}

\expandafter\ifx\csname >amsstdmathmacs\endcsname\relax
\expandafter\def\csname >amsstdmathmacs\endcsname{done}
\else\endinput\fi
\def\mathcs{C^{\displaystyle *}}
\def\cs{\ifmmode\mathcs\else$\mathcs$\fi}

\let\bbb=\Bbb
\def\R{{\bbb R}}
\def\C{{\bbb C}}
\def\T{{\bbb T}}
\def\Z{{\bbb Z}}
\def\K{{\cal K}}
\def\iff{if and only if}

\def\Ad{\operatorname{Ad}}

\def\Inn{\operatorname{Inn}}

\def\ker{\operatorname{ker}}
\def\tr{\operatorname{tr}}
\def\Aut{\operatorname{Aut}}
\def\sp{\operatornamewithlimits{span}}
\def\nbhd{neighborhood}
\def\id{\operatorname{id}}
\def\set#1{\{\,#1\,\}}
\let\tensor=\otimes
\def\({\bigl(}
\def\){\bigr)}
\def\restr#1{|_{{#1}}}
\def\spec#1{\specnp{(#1)}}
\def\specnp#1{{#1}^\wedge}
\newbox\hidebox
\def\spechide#1{\setbox\hidebox=\hbox{$#1$}
\hbox to\wd\hidebox{$\box\hidebox^\wedge$\hss}}
\expandafter\ifx\csname>ckrwmacs\endcsname\relax
\expandafter\def\csname>ckrwmacs\endcsname{done}
\else\endinput\fi
\expandafter\ifx\csname>amsspfonts\endcsname\relax
\expandafter\def\csname>amsspfonts\endcsname{done}
\else\endinput\fi
\makeatletter

\new@fontshape{spfonts}{m}{n}{%
   <5>1spfonts%
   <6>1spfonts%
   <7>1spfonts%
   <8>1spfonts%
   <9>1spfonts%
   <10>spfonts%
   <11>spfonts at10.95pt%
   <12>spfonts at12pt%
   <14>1spfonts at12pt%
   <17>1spfonts at12pt%
   <20>1spfonts at12pt%
   <25>1spfonts at12pt}{}

\new@fontshape{eus}{m}{n}{%
   <5>eusm5%
   <6>eusm6%
   <7>eusm7%
   <8>eusm8%
   <9>eusm9%
   <10>eusm10%
   <11>eusm10 at10.95pt%
   <12>eusm10 at12pt%
   <14>eusm10 at14.4pt%
   <17>eusm10 at17.28pt%
   <20>eusm10 at20.74pt%
   <25>eusm10 at24.88pt}{}

\extra@def{eus}{\skewchar#1'60}{}

\newmathalphabet*{\script}{eus}{m}{n}
\makeatother
\def\bimodfont#1{\script{#1}}
\def\Im{\operatorname{Im}}
\def\Homeo{\operatorname{Homeo}}
\def\Out{\operatorname{Out}}
\def\homeot{\Homeo(T)}
\def\X{{\bimodfont X}}
\def\Y{\bimodfont Y}

\def\smeover #1{\,\mathord{\mathop{\text{--}}\nolimits_{#1}}\,}
\def\sme{\,\mathord{\mathop{\text{--}}\nolimits_{\relax}}\,}
\def\eb{im\-prim\-i\-tiv\-ity bi\-mod\-u\-le}
\let\ipscriptstyle=\scriptscriptstyle
\def\lipsqueeze{{\mskip -3.0mu}}
\def\ripsqueeze{{\mskip -3.0mu}}
\newbox\ipstrutbox
\setbox\ipstrutbox=\hbox{\vrule height8.5pt% depth 3.5pt
width 0pt}
\def\ipstrut{\copy\ipstrutbox}
\def\lip#1<#2,#3>{\mathopen{\relax_{\ipstrut\ipscriptstyle{
#1}}\lipsqueeze
\langle} #2,#3 \rangle}
\def\blip#1<#2,#3>{\mathopen{\relax_{\ipstrut
\ipscriptstyle{ #1}}\lipsqueeze\bigl\langle} #2,#3 \bigr\rangle}
\def\rip#1<#2,#3>{\langle #2,#3
\rangle_{\ripsqueeze\ipstrut\ipscriptstyle{#1}}}
\def\brip#1<#2,#3>{\bigl\langle #2,#3
\bigr\rangle_{\ripsqueeze\ipstrut\ipscriptstyle{#1}}}
\def\angsqueeze{\mskip -6mu}
\def\smangsqueeze{\mskip -3.7mu}
\def\trip#1<#2,#3>{\langle\smangsqueeze\langle #2,#3
\rangle\smangsqueeze\rangle_{\ripsqueeze\ipstrut\ipscriptstyle{#1}}}
\def\btrip#1<#2,#3>{\bigl\langle\angsqueeze\bigl\langle #2,#3
\bigr\rangle
\angsqueeze\bigr\rangle_{\ripsqueeze\ipstrut\ipscriptstyle{#1}}}
\def\tlip#1<#2,#3>{\mathopen{\relax_{\ipstrut\ipscriptstyle{
#1}}\lipsqueeze \langle\smangsqueeze\langle} #2,#3
\rangle\smangsqueeze\rangle}
\def\btlip#1<#2,#3>{\mathopen{\relax_{\ipstrut\ipscriptstyle{
#1}}\lipsqueeze
\bigl\langle\angsqueeze\bigl\langle} #2,#3
\bigr\rangle\angsqueeze\bigr\rangle}
\def\btensor{\mathrel{\widehat\tensor}}
\def\Coftensor#1{\mathrel{\tensor_{#1}}}
\def\Ttensor{\Coftensor{C_0(T)}}
\def\bCoftensor#1{\mathrel{\widehat{\tensor}_{#1}}}
\def\bTtensor{\bCoftensor{C_0(T)}}
\def\cH{\check H}
\def\sT{\script S}
\def\gij{g_{ij}}
\def\hij{h_{ij}}
\def\hjk{h_{jk}}
\def\hik{h_{ik}}
\def\bgij{{\bar g}_{ij}}
\def\gik{g_{ik}}
\def\gjk{g_{jk}}
\def\fij{{F_{ij}}}
\def\fijk{{F_{ijk}}}
\def\sbrg{{\frak B \frak r}_G}%  modified 31-8-93
\def\br{\operatorname{Br}}%  added 31-8-93
\def\brg{\br_G}%  modified 31-8-93
\def\brt{\br(T)}%  modified 31-8-93
\def\sbrgt{\sbrg(T)}
\def\brgt{\brg(T)}
\def\tx{\tilde x}
\def\ty{\tilde y}
\def\dualell#1{(#1)^\sim}
\def\B{B}
\def\bX{\overline{\X}}
\def\bA{\overline{A}}
\def\hA{\hat{A}}
\def\atba{A\Ttensor\bA}
\def\hs\|#1\|{\|#1\|_{\scriptscriptstyle{\text{HS}}}}
\def\L{{\cal L}}
\def\m{{\frak M}}
\def\n{{\frak N}}
\def\xittxib{{\X_i\bCoftensor{C(F_i)}\bX_i}}
\def\c{{\frak c}}
\def\l{\lambda}
\def\g{\gamma}
\def\Autc{\Aut_{C_0(T)}}
\def\Outc#1{\Autc(#1)/\Inn(#1)}
\def\M{\cal M}
\def\U{\cal U}
\def\UM{\cal U \M}
\def\UZM{\cal U\cal Z \M}
\def\tv{\tilde v}
\def\ttau{\tau}% Changed 30/8/93
\def\w{\omega}
\def\blnorm{\Bigl\|}
\def\brnorm{\Bigr\|}
\theoremstyle{plain} %% This is the default
\newtheorem{spthm}{Theorem}

\newsymbol\rtimes 226F

\newif\ifaddresses
\addressesfalse

\begin{document}

\title[The Equivariant Brauer Group]{An Equivariant
Brauer Group and\\ Actions of Groups on \boldmath\cs-algebras}

%\dedicatory{This a prelimary version of this paper.  It should
%not be copied, circulated, or referenced.  The authors will be
%happy to supply a final version upon request.}

\ifaddresses
\else
\author[Crocker]{David Crocker}
\address{School of Mathematics \\
University of New South Wales \\
P.~O.~Box~1 \\
Kensington, NSW 2033 \\ Australia}

\author[Kumjian]{Alexander Kumjian}
\address{Department of Mathematics \\
University of Nevada \\
Reno, NV 89557 \\
USA}
\email{alex@@math.unr.edu}

\author[Raeburn]{Iain Raeburn}
\address{Department of Mathematics \\
University of Newcastle \\
Newcastle, NSW 2308 \\
Australia}
\email{iain@@frey.newcastle.edu.au}

\author[Williams]{\\Dana P. Williams}
\address{Department of Mathematics \\
Dartmouth College \\
Hanover, NH 03755-3551 \\
USA}
\email{dana.williams@@dartmouth.edu}
\fi

\thanks{The fourth author was partially supported by the
National Science Foundation.}

\thanks{This research was supported by the Australian
Research Council.}

\ifaddresses
\thanks{Please direct correspondence to the fourth author.}
\fi

\date{20 October 1993}



\maketitle

\ifaddresses
\begin{center}
\uppercase{David Crocker} \\
School of Mathematics,
University of New South Wales \\
P.~O.~Box~1,
Kensington, NSW 2033,
Australia \\[.5cm]
\uppercase{Alexander Kumjian}\\
Department of Mathematics,
University of Nevada\\
Reno, NV 89557,
USA \\
email: alex@@math.unr.edu \\[.5cm]
\uppercase{Iain Raeburn}\\
Department of Mathematics,
University of Newcastle\\
Newcastle, NSW 2308,
Australia \\
email: iain@@frey.newcastle.edu.au \\[.5cm]
\uppercase{Dana P. Williams}\\
Dartmouth College,
Department of Mathematics\\
Bradley Hall,
Hanover, NH 03755-3551,
USA\\
e-mail: dana.williams@@dartmouth.edu
\end{center}

\begin{abstract}
%\input ckrwab
%
%   SCCS Info: ckrwab.tex version 3.1: 10/19/93
%
Suppose that $(G,T)$ is a second countable locally compact
transformation group given by a homomorphism
$\ell:G\to\Homeo(T)$, and that $A$ is a
separable continuous-trace \cs-algebra with spectrum $T$.  An
action $\alpha:G\to\Aut(A)$ is said to cover $\ell$ if the
induced action of $G$ on $T$ coincides with the original one.
We prove that the set $\brgt$ of Morita
equivalence classes of such systems forms a group with
multiplication given by the
balanced tensor product: $[A,\alpha][B,\beta] = [A\Ttensor
B,\alpha\tensor\beta]$, and we refer to $\brgt$ as the
Equivariant Brauer Group.

We give a detailed analysis of the structure of $\brgt$ in
terms of the Moore cohomology of the group $G$ and the
integral cohomology of the space $T$. Using this, we can
characterize the stable continuous-trace \cs-algebras with
spectrum $T$ which admit actions covering $\ell$.  In
particular, we prove that if $G=\R$, then every stable
continuous-trace \cs-algebra admits an (essentially unique)
action covering~$\ell$, thereby substantially improving
results of Raeburn and Rosenberg.
\end{abstract}

\else

\begin{abstract}
%\input ckrwab
%
%   SCCS Info: ckrwab.tex version 3.1: 10/19/93
%
Suppose that $(G,T)$ is a second countable locally compact
transformation group given by a homomorphism
$\ell:G\to\Homeo(T)$, and that $A$ is a
separable continuous-trace \cs-algebra with spectrum $T$.  An
action $\alpha:G\to\Aut(A)$ is said to cover $\ell$ if the
induced action of $G$ on $T$ coincides with the original one.
We prove that the set $\brgt$ of Morita
equivalence classes of such systems forms a group with
multiplication given by the
balanced tensor product: $[A,\alpha][B,\beta] = [A\Ttensor
B,\alpha\tensor\beta]$, and we refer to $\brgt$ as the
Equivariant Brauer Group.

We give a detailed analysis of the structure of $\brgt$ in
terms of the Moore cohomology of the group $G$ and the
integral cohomology of the space $T$. Using this, we can
characterize the stable continuous-trace \cs-algebras with
spectrum $T$ which admit actions covering $\ell$.  In
particular, we prove that if $G=\R$, then every stable
continuous-trace \cs-algebra admits an (essentially unique)
action covering~$\ell$, thereby substantially improving
results of Raeburn and Rosenberg.
\end{abstract}

\fi

%
%   SCCS Info: ckrw1.tex version 3.1: 10/19/93
%
\section{Introduction}\label{intro}

In 1963, Dixmier and Douady associated to each continuous-trace
\cs-algebra $A$ with spectrum $T$ a class $\delta(A)$ in the
cohomology group $H^3(T;\Z)$, which determines $A$ up to a natural
equivalence relation \cite{dd,dixfields}. Over the past 15 years, it
has become clear that this relation is precisely the \cs-algebraic
version of Morita equivalence developed by Rieffel; this observation
appears, for example, in \cite{green4,beer}, and a modern treatment
of the theory is discussed in \cite[\S3]{90c}. It was also realized
in the mid--1970's that the results of \cite{dd,dixfields}
effectively establish an isomorphism between a Brauer group $\brt$
and $H^3(T;\Z)$: the elements of $\brt$ are Morita equivalence
classes $[A]$ of continuous-trace algebras $A$ with spectrum $T$, the
multiplication is given by the balanced \cs-algebraic tensor product
$[A][B]=[A\otimes_{C(T)}B]$, the identity is $[C_0(T)]$, and the
inverse of $[A]$ is represented by the conjugate algebra $\overline
A$. This point of view was discussed by Taylor \cite{joe} and Green
\cite{green4}, although neither published details.

Much of the current interest in operator algebras focuses on
\cs-dynamical systems, in which a locally compact group acts on
a \cs-algebra, and it is natural to try to extend the Dixmier-Douady
theory to accommodate group actions. Thus one starts with an action
of a locally compact group $G$ on a locally compact space $T$,
and considers systems $(A,\alpha)$ in which $A$ is a continuous-trace
\cs-algebra with spectrum $T$ and $\alpha$ is an action of $G$ on $A$
which induces the given action of $G$ on $T$. There is a notion of
Morita equivalence for dynamical systems due to Combes \cite{combes}
and Curto--Muhly--Williams \cite{cmw}, which is easily modified to
respect the identifications of spectra with $T$, and the elements of
our equivariant Brauer group $\brgt$ are the Morita equivalence
classes $[A,\alpha]$ of the systems $(A,\alpha)$. The group
operation is given by
$[A,\alpha][B,\beta]=[A\otimes_{C(T)}B,\alpha\otimes_{C(T)}\beta]$,
the identity is $[C_0(T),\tau]$, where $\tau_s(f)(x)=f(s^{-1}\cdot
x)$, and the inverse of $[A,\alpha]$ is $[\overline
A,\overline\alpha]$, where $\overline\alpha(\overline
a):=\overline{\alpha(a)}$. Even though the key ideas are all in
\cite{dixfields}, it is not completely routine that $\brgt$ is a
group, and we have to work quite hard to establish that
$(A\otimes_{C(T)}\bA,\alpha\otimes_{C(T)}\overline\alpha)$
is Morita equivalent to $(C_0(T),\tau)$.

Similar Brauer groups have been constructed by Parker
for $G=\Z/2\Z$ \cite{emp},
and by
Kumjian in the context of $r$-discrete
groupoids
\cite{alex3}.  The results of the preceding paragraph are contained
in those of
\cite{alex3} when the group is discrete. However, Kumjian then
generalizes the Dixmier--Douady Theorem by identifying his Brauer
group with the equivariant cohomology group $H^2(T,G;\sT)$ of
Grothendieck \cite{groth}. (If $G$ is trivial, $H^2(T,\sT)$ is
naturally isomorphic to $H^3(T;\Z)$, and the original Dixmier--Douady
construction proceeds through $H^2(T,\sT)$.)  Grothendieck developed
powerful techniques for computing his equivariant cohomology, and
there is in particular a spectral sequence $\{E^{p,q}_r\}$ with
$E^{p,q}_2=H^p\(G,H^q(T,\sT)\)$ (the group cohomology of $G$ with
coefficients in the sheaf cohomology of $T$) which converges to
$H^{p+q}(T,G;\sT)$. In view of Kumjian's result, this gives a
filtration of the equivariant Brauer group $\brgt$ for discrete $G$.

For locally compact groups, the appropriate version of group
cohomology is the Borel cochain theory developed by Moore
\cite{moore3,moore4}. (Computing the 2-cocycle for the extension
$0\to\Z\to\R\to\T\to1$ shows that continuous cochains will not
suffice.) The coefficient modules in Moore's theory must be Polish
groups, and there are not enough injective objects in this category
to allow the direct application of homological algebra, so any
suitable generalization of Grothendieck's theory will, at best, be
hard to work with. However, we are only interested here in the
Brauer group, and the filtration involves only the low-dimensional
cohomology groups $H^p\(G,H^q(T,\sT)\)$ for $p=0,1,2,3$ and
$q=0,1,2$. Each of the coefficient groups $H^0(T,\sT)=C(T,\T)$,
$H^1(T,\sT)\cong H^2(T;\Z)$ and $H^2(T,\sT)\cong H^3(T;\Z)$ admits a
\cs-algebraic interpretation: $H^2(T,\sT)$ is itself the Brauer
group of continuous-trace algebras with spectrum $T$, $H^1(T,\sT)$ is
the group $\operatorname{Aut}_{C_0(T)}A/\operatorname{Inn}A$ of outer
$C(T)$-automorphisms of a stable continuous-trace algebra $A$ with
spectrum $T$ \cite{pr1}, and $C(T,\T)$ is the unitary group $UZM(A)$
of the center of the multiplier algebra $M(A)$ of such an algebra $A$.
Further, the Moore cohomology groups $H^2\(G,C(T,\T)\)$ and
$H^3\(G,C(T,\T)\)$ arise naturally in the analysis of group actions on
a continuous-trace algebra $A$ with spectrum $T$: $H^2$ contains the
obstructions to implementing an action $\alpha:G\to \Inn (A)$ by a
unitary group
$u:G\to UM(A)$ \cite[\S0]{rr}, and $H^3$ the obstructions to
implementing a homomorphism $\beta:G\to \Aut (A)/ \Inn (A)$ by
a twisted action (see \lemref{lem:5.5} below).
The remarkable point of the present paper
is that, using these interpretations, we have been able to define all
the groups and homomorphisms necessary to completely describe the
filtration of $\brgt$ predicted by the isomorphism $\brgt\cong
H^2(T,G,\sT)$ of the discrete case. Thus we will prove:

\begin{spthm}[cf. \thmref{thm:structure} below]
Let $(G,T)$ be a second countable locally compact transformation
group with $H^2(T;\Z)$ countable. Then there is a composition series
$\{0\}\leq B_1\leq B_2\leq B_3=\brgt$ of the equivariant Brauer
group in which $B_3/B_2$ is isomorphic to a subgroup of $H^3(T;\Z)$,
$B_2/B_1$ to a subgroup of $H^1\(G,H^2(T;\Z)\)$, and $B_1$ to a
quotient of $H^2\(G, C(T,\T)\)$. Further, we can precisely identify the
subgroups and quotients in terms of homomorphisms between groups of
the form $H^p\(G,H^q(T;\Z)\)$.
\end{spthm}

The sting of this theorem lies in, first, the specific nature of the
isomorphisms, and, second, in the last sentence, where the
homomorphisms are all naturally defined using the \cs-algebraic
interpretations of Moore and \v Cech cohomology. The isomorphism $F$
of $B_3/B_2$ into $H^3(T;\Z)$ takes $[A,\alpha]$ to $\delta(A)$, so
its kernel $B_2$ is the set of classes of the form
$[C_0(T,\K),\alpha]$. For the isomorphism of $B_2$ into
$H^1\(G,H^2(T;\Z)\)$, we use the exact sequence
\begin{ediagram*}\dgARROWLENGTH=1em
\node{0}\arrow{e}\node{\Inn{C(T,\K)}}\arrow[2]{e}
\node[2]{\Aut_{C(T)}{C(T,\K)}}\arrow[2]{e,t}{\zeta}\node[2]{H^2(T;\Z)}
\arrow{e}\node{0}
\end{ediagram*}
of \cite{pr1}, and send $(C_0(T,\K),\alpha)\in B_2$ to the cocycle
$s\mapsto\zeta(\alpha_s)$ in $Z^1\(G,H^2(T;\Z)\)$. Thus $B_1$ consists
of the systems $(C_0(T,\K),\alpha)$ in which $\alpha:G\to
\Inn{C_0(T,\K)}$, and the last isomorphism takes such an action
$\alpha$ to its Mackey obstruction---the class in
$H^2\(G,H^2(T;\Z)\)$ which vanishes precisely when $\alpha$ is
implemented by a unitary group $u:G\to UM(A)$.

To illustrate the second point, we describe our identification of the
range of the first homomorphism $F:\brgt\to H^3(T;\Z)$. We first
restrict attention to the group $H^3(T;\Z)^G$ of classes fixed under
the canonical action of $G$, and define a homomorphism
$d_2:H^3(T;\Z)^G\to H^2\(G,H^2(T;\Z)\)$. We then define another
homomorphism $d_3$ from the kernel of $d_2$ to a quotient of
$H^3\(G,C(T,\T)\)$, and prove that the image of $F$ is the kernel of
$d_3$. To see why this is powerful, note that a stable algebra $A$
with spectrum $T$ carries  an action of $G$ covering the given action
on $T$ if and only if $\delta(A)\in\Im F$. When $G=\R$,
$H^3(T;\Z)^\R=H^3(T;\Z)$, and results from \cite{rr} show that
$H^3\(\R,C(T,\T)\)=H^2\(\R,H^2(T;\Z)\)=0$; we deduce that $F$ maps onto
$H^3(T;\Z)$, and hence that {\it every} action of $\R$ on $T$ lifts
to an action of $\R$ on {\it every} stable  continuous-trace algebra
$A$ with spectrum $T$ (see Corollary~\ref{cor:5.9} below).
This is a substantial generalization of results proved
in \cite[\S4]{rr}---and even they required considerable machinery.

We should stress that, even when there is no group action and $T$ is
compact, our Brauer group $\brt$ is not the usual Brauer group of the
commutative ring $C(T)$, which is isomorphic to the torsion subgroup of
$H^3(T;\Z)$ rather than $H^3(T;\Z)$ \cite{groth2}. Although the two
groups have different objects, $\brt$ is
isomorphic to the bigger Brauer group $\widetilde B\(C(T)\)$ of Taylor
\cite{joe2,rt2}, which is a purely algebraic invariant designed to
accommodate non-torsion cohomology classes. Presumably there is also an
equivariant version of $\widetilde B(R)$ for which theorems similar to
ours are true---indeed, the results in \cite{rt2,alex3} suggest that
$\widetilde B_G(R)$ might then be isomorphic to an equivariant \'etale
cohomology group.

\bigbreak
Our work is organized as follows.  In \secref{sec:2} we outline
some of the basic definitions of the internal and external tensor
products of \eb s which are fundamental to our approach.  In
\secref{sec:3} we discuss the Morita equivalence of systems, define
our Brauer group, and prove that it is indeed a group. We then devote
\secref{sec:5} to identifying the range of our Forgetful Homomorphism
$F:\brgt\to \brt\cong H^3(T;\Z)$, which is probably the most
important part of our main theorem. In \secref{sec:structure},
we give
a precise statement of our theorem, and finish off its proof. In the
last section, we discuss the application to actions of $\R$, and
consider some special cases in which we can say more about $\brgt$.

\bigbreak
We will adopt the following conventions.  When we consider a
{\it \cs-algebra $A$ with spectrum $T$} we are
considering a pair $(A,\phi)$ where $\phi:\hA\to T$ is a fixed
homeomorphism.  While we have opted to be less pedantic and drop the
$\phi$, it is necessary to keep its existence
in mind.  Thus, as in \cite[\S2]{90c}, we will work almost
exclusively with complete \eb s which preserve
the spectrum: if $A$ and $B$ are \cs-algebras with
spectrum $T$, then $\X$ is an $A\smeover T B$-\eb{}
if $\X$ is an \eb{} in the usual sense and the Rieffel
homeomorphism $h_\X:T\to T$ is the identity.  It is convenient to
keep in mind that, if $A$ and $B$ have continuous trace, then
if follows from Proposition~1.11 and the preceding remarks in
\cite{ra1} that $h_\X=\id$ if and only if the left and right actions
of $C_0(T)$ on $\X$, induced by the actions of $A$ and
$B$, respectively, coincide: i.e., $\phi\cdot x =x\cdot \phi$ for all
$\phi\in C_0(T)$ and $x\in\X$.
 (See \cite[\S2]{90c} for further
details.) We will also make full use of dual \eb s as defined in
\cite[Definition~6.17]{rieff1}.  Recall that if $\X$ is an $A\smeover
TB$-\eb, the dual $\widetilde\X$ of $\X$ is the set $\{\tilde
x:x\in\X\}$, made into a $B\smeover TA$-\eb\ as follows:
\begin{alignat*}{2} b\cdot\tx&=\dualell{x\cdot b^*}&\qquad \tx\cdot
a&=\dualell{a^*\cdot x} \\
\lip B<\tx,\ty>&=\rip B<x,y>  &\qquad
\rip A<\tx,\ty>&=\lip A<x,y>,
\end{alignat*}
for $x,y\in\X$, $a\in A$, and $b\in B$.

We will use the notation $H^n(T;\Z)$ for the ordinary integral
cohomology groups, and $H^n(T,\sT)$ for the sheaf cohomology
groups with coefficients in the sheaf of germs of
continuous circle-valued functions on $T$.  We will make frequent use
of the canonical isomorphism of $H^2(T,\sT)$ and $H^3(T;\Z)$; in
particular, we will view the Dixmier-Douady class $\delta(A)$ of a
continuous-trace \cs-algebra with spectrum $T$ as belonging to
whichever of these groups is more convenient for the matter at
hand.  It will also be essential to use Moore's Borel cochain
version of group cohomology as presented in \cite{moore3}: when $G$ is
a locally compact group, and $A$ is a Polish $G$-module,
$H^n(G,A)$ will denote the corresponding Moore group.

\bigbreak
The construction of our equivariant Brauer group was originally
intended to be part of the first author's Ph.D. thesis; in particular,
\thmref{thm:crocker} is basically due to
him.
Much of this work was carried out while the first three
authors were at the University of New South Wales. It was finished
while the third author was visiting the University of Colorado, and he
thanks his colleagues there for their warm hospitality. This research
was supported by the Australian Research Council.

%
%   SCCS Info: ckrw2.tex version 3.1: 10/19/93
%
\section{Tensor products of \eb s}
\label{sec:2}

Let $A$, $B$, $C$, and $D$ be \cs-algebras.
Suppose that $\X$ is a $A\sme B$-\eb{} and that $\Y$ is a
$B\sme C$-\eb.  Then the algebraic tensor product $\X\odot\Y$
is a $A\sme C$-bimodule and carries $A$- and $C$-valued inner
products defined, respectively, by
\begin{gather}
\trip C<x\otimes y,x'\otimes y'>=\brip C<{\rip B<x',x>}y,y'>\\
\tlip A<x\otimes y,x'\otimes y'>= \blip A<x,x'{\lip
B<y',y>}>.
\end{gather}
It is straightforward to verify that $\X\odot\Y$ is a (pre-) $A\sme
C$-\eb, and
we shall write $\X\otimes_B\Y$ for the completion
with respect to the common semi-norm induced by the inner
products (see \cite[\S3]{rieff2}).
Suppose that in addition $A$, $B$, and $C$ have spectrum
(identified with) $T$, and that
$\X$ is a $A\smeover T
B$-\eb{} and $\Y$ is a $B\smeover T C$-\eb.
Then it
is shown in \cite[Lemma~1.3]{ra1}, that
$\X\tensor_B \Y$
is a $A\smeover T C$-\eb.  Although $\X\tensor_B\Y$ is not a
Banach space tensor product in the usual sense, it does follow
from \cite[Proposition~2.9]{rieff1} that
\begin{equation}\label{eq:bics}
\|x\tensor y\|\le\|x\|\,\|y\|.
\end{equation}

The construction above is an example of an {\it internal tensor
product\/} of Hilbert modules as described in \cite[\S1.2]{je-th}.
We will also need the {\it external tensor product}.
Specifically, if $\X$ is a $A\sme C$-\eb{} and $\Y$ is a $B\sme
D$-\eb, then the formulas
\begin{gather}
\tlip A\tensor B<x\tensor y,x'\tensor y'> =
\lip A <x, x'>\tensor \lip B<y,y'> \label{eq:btensor1}\\
\trip C\tensor D<x\tensor y,x'\tensor y'> =
\rip C <x, x'>\tensor \rip D<y,y'>\label{eq:btensor2}
\end{gather}
define, respectively, $A\odot B$- and $C\odot D$-valued
sesqui-linear forms
 on $\X\odot\Y$.  It follows from \cite[1.2.4]{je-th} that these
forms are inner products for any \cs-norms on $A\odot B$ and
$C\odot D$, and that in particular $\X\odot \Y$ can be
completed
to a $A\tensor B \sme  C\tensor D$-\eb{} (recall that `$\tensor$'
denotes that minimal tensor product)\footnote{This
is observed in \cite[\S13.5]{black}, and in
\cite[Proposition~2.9]{bui-thesis} where it is also observed that
the same holds
for the maximal tensor product.  In general, if $\nu$
is a \cs-norm on $C\odot D$, then $A\odot B$ acts
as adjointable bounded operators with respect to the right Hilbert
$C\tensor_\alpha D$-module
structure on $\X\odot\Y$.  This provides $A\odot B$
with a \cs-norm $\nu^*$ for which the completion of $\X\odot \Y$ is
a $A\tensor_{\nu^*}B\sme C\tensor_\nu D$-\eb.  Since all our
algebras  will be
continuous-trace \cs-algebras, and hence nuclear,  the result from
\cite{je-th} will suffice.}. In order to more clearly distinguish
which tensor product of \eb s we're using,  we shall  write
$\X\btensor\Y$ for the completion of $\X\odot\Y$ with respect to
the operations in
\eqref{eq:btensor1}~and
\eqref{eq:btensor2}.

Now suppose that $A$, $B$, $C$, and $D$ have
Hausdorff
spectrum $T$ and
that $\X$ is a $A\smeover T C$-\eb{} and $\Y$ a $B\smeover T D$-\eb.
In particular, by the Dauns-Hofmann
 Theorem,
$C_0(T)$ sits in the center of the multiplier
algebras of all these algebras so that $\X$ and $\Y$ are
$C_0(T)$-bimodules.  Therefore we can form the balanced tensor
products $A\Ttensor B$ and $C\Ttensor D$.  Each of these algebras
has spectrum $T$ (cf., e.g., \cite[Lemma~1.1]{rw}).
Recall that $A\Ttensor B$ is the quotient of $A\tensor B$ by the
closed
ideal $I_T$ spanned by $\set{\phi\cdot a\tensor b-a \tensor
\phi\cdot b:\text{$a\in A$, $b\in B$, and, $\phi \in C_0(T)$}}$.
Similarly, $C\Ttensor D$ is the quotient of $C\tensor D$ by an
ideal $J_T$.
\begin{lem}
Suppose that $A$, $B$, $C$, and $D$ are \cs-algebras with
Hausdorff spectrum $T$, and that $\X$ and $\Y$ are, respectively,
$A\smeover T B$- and $C\smeover T D$-\eb s.  Then the
correspondence of \cite[\S3]{rieff2} between ideals in $C\tensor D$
and ideals in $A
\tensor B$, induced by $\X\btensor\Y$,
maps $I_T$ to $J_T$.  In particular, the corresponding quotient $
\X\bTtensor\Y$ of $\X\btensor\Y$ is a $A\Ttensor B\smeover T
C\Ttensor D$-\eb.
\end{lem}
\begin{pf}
Let $K(I_T)$ be the ideal of $C\tensor D$ corresponding to
$I_T$ via the Rieffel correspondence $K$.
Since $I_T\cdot(\X\btensor\Y)$ is spanned by vectors of the form
$$(\phi\cdot a\tensor b - a\tensor \phi\cdot b)\cdot(x\tensor y)=
\((\phi\cdot a)\cdot x \tensor b\cdot y\) - \(a\cdot x\tensor ( \phi
\cdot b)\cdot y\),$$
where $a\in A$, $b\in B$, $\phi\in C_0(T)$, $x\in \X$, and $y\in
\Y$, it follows that
$$K(I_T)\subseteq \overline{\sp}\set{\trip C\tensor D<v,u>: v\in V_0,
u\in \X\odot\Y},$$
where $V_0=\sp\set{\phi\cdot x\tensor y - x\tensor \phi\cdot y:
\phi\in C_0(T), x\in\X, y\in\Y}$.
Consequently, $K(I_T)\subseteq J_T$.  By symmetry,
we have $J_T\subseteq
K(I_T)$, and therefore $K(I_T)=J_T$, which proves the first
assertion.

The second assertion will follow from the first and the
discussion preceding the lemma once we show that the left
and right
$C_0(T)$-actions on the quotient module
$(\X\btensor\Y)/I_T\cdot(\X\btensor\Y)$ coincide.
But
$$\phi\cdot[x\tensor y]=[\phi\cdot x\tensor y]=
[x\cdot\phi\tensor y]
=[x\tensor y]\cdot \phi.$$
(The first equality holds because
$\phi\cdot(a\tensor b)=(\phi\cdot a\tensor b)=(a\tensor \phi\cdot
b)$ in $A\Ttensor B$, and the second because the
module
$\X$ is
$T$-balanced by assumption.  The third is similar to the
first.)
\end{pf}

The next result is implicit in \cite{dixfields}.
Our approach here views the Dixmier-Douady class $\delta(A)$
of a continuous-trace \cs-algebra $A$ as the obstruction to
the existence of a global Morita equivalence of $A$ with $C_0(T)$
as described in \cite[\S3]{90c}.
\begin{prop}\label{prop:2.2}
Suppose that $A$ and $B$ are continuous-trace \cs-algebras with
spectrum $T$.  Then $\delta(A\Ttensor B)=\delta(A)+\delta(B)$.
\end{prop}
\begin{pf}
Since $A$ has continuous trace, it follows
from \cite[Lemmas 6.1~and 6.2]{90c}
that there are compact sets $F_i\subseteq
T$ whose interiors form a cover $\frak
A=\set{\operatorname{int}F_i:i\in I}$ of $T$  such that:
\begin{enumerate}
\item
for each $i\in I$ there are $A^{F_i}\smeover{F_i}C(F_i)$-\eb s
$\X_i$, and
\item
for each $i,j\in I$, there are \eb{} isomorphisms $\gij:\X_j^\fij\to
\X_i^\fij$.
\end{enumerate}
Then the class $\delta(A)$  in $H^3(T,\Z)$ is determined by the
cocycle $\nu=\set{\nu_{ijk}}$ in $\cH^2({\frak A},\sT)$ defined
by
$$\gij^\fijk\(\gjk^\fijk(x)\)=\nu_{ijk}\cdot \gik^\fijk(x).$$
By taking refinements, we may assume that we have
similar data for $B$ consisting of bimodules $\set{\Y_i}$,
isomorphisms $\set{h_{ij}}$, and a cocycle $\mu=\set{\mu_{ijk}}$
all defined with respect to the same cover $\frak A$.

The result follows from verifying that $(A\Ttensor B)^{F_i}\cong
A^{F_i}\Coftensor{C(F_i)} B^{F_i}$, that
$(\X_i\bCoftensor{C(F_i)}\Y_i)^\fij\cong
\X_i^\fij\bCoftensor{C(\fij)}\Y_i^\fij$, and that $k_{ij}=\gij\tensor
h_{ij}$ defines an isomorphism of
$\X_j^\fij\bCoftensor{C(\fij)}\Y_j^\fij$ onto
$\X_i^\fij\bCoftensor{C(\fij)}\Y_i^\fij$.
Then  $k_{ij}^\fijk\circ k_{jk}^\fijk = \nu_{ijk}\mu_{ijk}\cdot
h_{ik}^\fijk$.
\end{pf}
%
%   SCCS Info: ckrw3.tex version 3.1: 10/19/93
%
\section{The Brauer Group}
\label{sec:3}

For
the remainder of this article, $(G,T)$ will be a second
countable locally compact transformation group. We define
$\sbrgt$ to be the class of pairs $(A,\alpha)$ where $A$ is a
separable continuous-trace \cs-algebra with spectrum $T$ and
$\alpha:G\to\Aut(A)$ is a strongly continuous action inducing the
given action $\tau$ on $C_0(T)$.  That is, $\alpha_s(\phi\cdot
a)=\tau_s(\phi)\cdot \alpha_s(a)$ for $a\in A$ and $\phi\in
C_0(T)$, where $\tau_s(\phi)(t)=\phi(s^{-1}\cdot t)$.

We say that two elements $(A,\alpha)$ and $(B,\beta)$ of $\sbrgt$
are equivalent, written $(A,\alpha)\sim(B,\beta)$, if they are
Morita equivalent over $T$
in the sense of Combes \cite{combes} (see also \cite[\S4]{90c}):
this means that there is an
$A\smeover T B$-\eb{} $\X$ and an action $u$ of {$G$} on $\X$ by linear
transformations,  which is strongly
continuous (i.e., $s\mapsto u_s(x)$ is norm-continuous for all
$x$) and
satisfies
\begin{equation*}
\alpha_s\(\lip A<x,y>\)=\blip A<u_s(x),u_s(y)>\quad
\text{and}\quad
\beta_s\(\rip B<x,y>\)=\brip
B<u_s(x),u_s(y)>.
\end{equation*}

We claim that $\sim$ is an equivalence relation.  It is certainly
reflexive: take $(\X,u)=(A,\alpha)$. Symmetry is immediate
from the existence of dual \eb s: one only has to define $\tilde
u_s$ by $\tilde u_s(\tx)=\dualell{u_s(x)}$.
Transitivity requires more work.  Suppose that
$(A,\alpha)\sim(B,\beta)$ via $(\X,u)$ and that $(B,\beta) \sim
(C,\gamma)$ via $(\Y,v)$.  Then we have
\begin{align}\blnorm u_s\tensor v_s\Bigl(\sum_i x_i\tensor
y_i\Bigr)\brnorm^2 &=
\blnorm\sum\nolimits_i u_s(x_i)\tensor v_s(y_i)\brnorm^2
\notag\\ &=
\blnorm\sum\nolimits_{i,j}\brip C<{\brip B <u_s(x_j),u_s(x_i)>}
v_s(y_i),v_s(y_j)>\brnorm \notag\\
&=\blnorm\sum\nolimits_{i,j}\brip C<\beta_s\({\rip B<x_j,x_i>}\)\cdot
 v_s(y_i),
v_s(y_j)>\brnorm \label{eq:w}\\
&=\blnorm\sum\nolimits_{i,j}\brip C<v_s\({\rip B<x_j,x_i>}\cdot y_i\),
 v_s(y_j)>
\brnorm \notag\\
\intertext{which, because $v_s$ is isometric, is}
&=\blnorm\sum\nolimits_{i,j}\brip C<{\rip B<x_j,x_i>}\cdot y_i,
y_j>\brnorm
\notag\\
&= \blnorm\sum\nolimits_i x_i\tensor y_i\brnorm^2 .\notag
\end{align}
Therefore $w_s=u_s\tensor v_s$ defines an
action of $G$ on $\X\tensor_B\Y$, which is strongly
continuous in view of \eqref{eq:bics},
and $(\X\tensor_B\Y,w)$ provides the required
equivalence between $(A,\alpha)$ and $(C,\gamma)$.
We will write $\brgt$ for the set $\sbrgt/\!\!\sim$
of equivalence classes\footnote{Notice that $\brgt$ is actually a
set.  To see this, fix an infinite dimensional separable Hilbert space
$\cal H$.  Then the separable \cs-subalgebras of $\B(\cal H)$ form a
set as do all possible actions of $G$ on each subalgebra, and
therefore as do the collection of all possible actions on separable
subalgebras of $\B(\cal H)$.  This set contains a representative
for each equivalence class.  (Separability could be
replaced by any fixed cardinality.)}.

It will be helpful to keep in mind that the above equivalence
relation can be reformulated as follows.
Recall that two actions $\alpha:G\to\Aut(A)$ and $\beta:G\to\Aut(B)$ are
{\em outer conjugate\/} if there is an isomorphism $\Phi:A\to B$ so that
$\beta$ is exterior equivalent to $\Phi\circ\beta\circ\Phi^{-1}$.
We say that $\alpha$ and $\beta$ are {\em stably outer conjugate\/} if
$\alpha\tensor i$ and $\beta\tensor i$ are outer conjugate as actions on
$A\tensor\K$ and $B\tensor\K$, respectively.
If $A$ and $B$ have spectrum $T$, then we say that $\alpha$ is outer
conjugate over $T$ if $\Phi$ can be taken to $C_0(T)$-linear.
\begin{lem}\label{lem:combes}
Suppose that $(A,\alpha),(B,\beta)\in\sbrgt$.  Then
$(A,\alpha)\sim(B,\beta)$ if and only if $\alpha$ is stably outer
conjugate to $\beta$ over $T$.  If $A$ and $B$ are both stable, then
$(A,\alpha)\sim(B,\beta)$ if and only if $\alpha$ is outer conjugate to
$\beta$ over $T$.
\end{lem}
\begin{pf}
The first statement follows from the second since
$(A,\alpha)\sim(A\tensor\K,\alpha\tensor i)$ and $(B,\beta)\sim
(B\tensor\K,\beta\tensor i)$.
But, if $A$ and $B$ are stable, and $(A,\alpha)\sim(B,\beta)$, then the
proposition in \S9 of \cite{combes} implies that $\alpha$ and $\beta$ are
outer conjugate.  The argument in the last paragraph of the proof of
Lemma~2.3 in \cite{90b} shows that the isomorphism of $A$ onto $B$ can be
taken to be $C_0(T)$-linear.  Finally, if $\alpha$ and $\beta$ are outer
conjugate over $T$, then the other half of the same proposition
implies that $(A,\alpha)$ is Morita equivalent to $(B,\beta)$, and
again, it is straightforward to see that we can take the Morita
equivalence over $T$.
\end{pf}

Let $(A,\alpha)$ and $(B,\beta)$ be elements of $\sbrgt$.
Notice that
$\alpha_s\tensor \beta_s(\phi\cdot a \tensor b - a\tensor \phi
\cdot b)= (\phi\circ\tau_s^{-1})\cdot\alpha_s(a)\tensor
\beta_s(b) - \alpha_s(a)\tensor (\phi\circ\tau_s^{-1})\cdot
\beta_s(b)$.  Thus, $\alpha_s\tensor\beta_s$ maps the closed
ideal $I_T $ of $A\tensor B$ spanned by $\set{
\phi\cdot a\tensor b-a\tensor \phi\cdot b}$
to
itself, and defines an automorphism $\alpha_s \Ttensor \beta_s$
of $A\Ttensor B=A\tensor B/I_T$.
It is easy to check that $\alpha_s
\Ttensor \beta_s$ induces the given action on $T$, so that
$(A\Ttensor B,\alpha\Ttensor \beta)\in\sbrgt$.  For notational
convenience, we will usually write $\alpha\tensor\beta$ rather
than $\alpha\Ttensor\beta$.
\begin{lem}\label{lem:wd}
Suppose that $(A,\alpha)\sim(C,\gamma)$ via $(\X,u)$ and that $
(B,\beta)\sim(D,\delta)$ via $(\Y,v)$.  Then $(A\Ttensor
B,\alpha\tensor\beta)$ is equivalent to $(C\Ttensor
D,\gamma\tensor\delta)$ in $\sbrgt$.
\end{lem}
\begin{pf}
As pointed out in \secref{sec:2}, $\X\bTtensor\Y$ is an
$A\Ttensor B \smeover T C\Ttensor D$-\eb.  The argument that $w_s(
x\btensor y)=u_s(x)\btensor v_s(y)$ gives a well-defined strongly
continuous action of $G$ on $\X\bTtensor\Y$ is similar, but more
straightforward, than \eqref{eq:w} above.
Then $(\X\bTtensor\Y,w)$ is the required $(A\Ttensor
B,\alpha\tensor\beta)\smeover T (C\Ttensor
D,\gamma\tensor\delta)$-\eb.
\end{pf}
\begin{prop}
The binary operation
\begin{equation}\label{eq:binop}
[A,\alpha][B,\beta]=[A\Ttensor B,\alpha\tensor\beta]
\end{equation}
is well-defined on $\brgt$, and with respect to this operation,
$\brgt$ is a commutative semi-group with identity equal to the class
of $(C_0(T),\tau)$.
\end{prop}
\begin{pf}
The operation \eqref{eq:binop} is well-defined by virtue of
\lemref{lem:wd}.  Since an equivariant
$C_0(T)$-isomorphism of $A$ onto $B$ certainly gives a
Morita equivalence over $T$, associativity and commutativity follow from
the observations that $(a\tensor b)\tensor c\mapsto a\tensor(b\tensor
c)$ and $a\tensor b\mapsto b\tensor a$ define equivariant
$C_0(T)$-isomorphisms of
$(A\Ttensor B)\Ttensor C$ onto $A \Ttensor(B\Ttensor C)$ and $A\Ttensor
B$ onto $B\Ttensor A$, respectively.  Similarly,
$(C_0(T),\tau)$ is an identity because the isomorphism $a\tensor
f\mapsto f\cdot a$ of
$A\Ttensor C_0(T)$ onto $A$ is equivariant and  $C_0(T)$-linear.
\end{pf}
\begin{remark}\label{rem:3.3}
If $G=\set e$, we write $\brt$ for $\brgt$.  It is well-known
that the map
sending $[A]$ to $\delta(A)$ defines a bijection of $\brt$
with $H^3(T;\Z)$ (see, for example, \cite[Theorem~3.5]{90c}).
\propref{prop:2.2} implies that
$[A]\mapsto\delta(A)$ is a (semi-group) isomorphism; in
particular, $\brt$ is a group.
\end{remark}
%
%   SCCS Info: ckrw4.tex version 3.1: 10/19/93
%
%\section{The Brauer Group}
%\label{sec:crocker}

If $V$ is a complex vector space, then we will write $\overline{V}
$ for the conjugate space: that is, $\overline{V} $ coincides
with $V$ as a set, and if $\flat=\flat_V:V\to\overline{V}$ is the
identity map, then scalar multiplication on $\overline{V}$ is
given by $\lambda\cdot\flat(v)=\flat(\bar\lambda\cdot b)$.
In the event $A$ is a \cs-algebra with spectrum $T$, then $\overline
{A}$ is again a \cs-algebra\footnote{In fact,
$\bA$ is $*$-isomorphic to the ``opposite algebra''
$A^{\text{op}}$ via the map $x^{\text{op}}\in A^{\text{op}}\mapsto
\flat(x^*)\in\bA$.}
with spectrum $T$, and if $\phi\in
C_0(T)$, then $\phi\cdot\flat_A(a)=\flat_A(\bar\phi\cdot a)$.
Furthermore, if $\X$ is an $A\smeover TB$-\eb, then $\bX$ is
naturally an $\overline A \smeover T \overline B$-\eb:
\begin{alignat*}{2}
\flat_A(a)\cdot\flat_\X(x)&=\flat_\X(a\cdot x)
&\qquad \flat_\X(x)\cdot\flat_B(b)&= \flat_\X(x\cdot b) \\
\blip \overline{A}<\flat_\X(x),\flat_\X(y)> & =\flat_A\( \lip A <x,y>\)
&\qquad \brip \overline B <\flat_\X(x),\flat_\X(y)> &=
\flat_B\(\rip B<x,y>\) .
\end{alignat*}
Of course, if $(A,\alpha)$ is in $\sbrgt$, then so is $(\overline
A,\bar\alpha)$, where $\bar\alpha_s\(\flat(a)\) = \flat\(\alpha_s
(a)\)$.
\begin{remark}
If $\X$ is a $A\smeover T C_0(T)$-\eb, then we will view $\X
\bTtensor\bX$ as a $A\Ttensor \bA \smeover T C_0(T)$-\eb{} by
identifying $C_0(T)\Ttensor \overline{C_0(T)}$ with $C_0(T)$ via the
isomorphism determined by $\phi\tensor\flat(\psi)\mapsto
\phi\overline{\psi}$.
\end{remark}
\begin{thm}\label{thm:crocker}
With the binary operation defined in \eqref{eq:binop}, $\brgt$ is
a group.  The inverse of $[A,\alpha]$ is given by $[\overline A,
\bar \alpha]$.
\end{thm}

\begin{remark}\label{rem:4.10}
The theorem has several immediate and interesting consequences.  For
example, we can reduce the problem of classifying $G$-actions on a
given stable  continuous-trace \cs-algebra $A$ with spectrum $T$
which cover the given action $\ell$ on $T$ to (1)~finding an {\em
single\/} action $\alpha$ on $A$ covering $\ell$
and (2)~classifying all $G$-actions on $C_0(T,\K)$ covering $\ell$.
To make this precise,  observe that the homomorphism
$F:\brgt\to\brt$ defined by $F\([A,\alpha]\)=\delta(A)$ (called the
{\it Forgetful  Homomorphism}) has as its kernel exactly the subgroup
of $\brgt$ consisting of classes (which have representatives) of the
form $\(C_0(T,K),\sigma\)$. Then the assertion above is simply that
the classes in $\brgt$ coming from
actions on $A$ are precisely
those in
$F^{-1}\(\delta(A)\)=[A,\alpha]\ker(F)$.
\end{remark}
\medbreak
To prove \thmref{thm:crocker}, all
that remains to be shown is the last assertion.  This will require the
remainder of the section.
We fix $(A,\alpha)$ in $\sbrgt$.  As before, we can choose data
$\set{F_i}_{i\in I}$, $\set{\X_i}$, $\set{\gij}$, and $\set{\nu_{ijk}}$
as in \propref{prop:2.2}.
Naturally, we can define $\bgij:\bX_J^\fij\to\bX_i^\fij$
by $\bgij\(\flat(x)\)=\flat\(\gij(x)\)$.  Then we can produce
data for $\atba$ of the form $\set{F_i}$, $\set{\xittxib}$,
and $\set{h_{ij}}$, where $h_{ij}=\gij\btensor \bgij$.  Notice
that
\begin{equation}\label{eq:hij}
\hij^\fijk\circ\hjk^\fijk=\hik^\fijk.
\end{equation}
Using the cocycle property \eqref{eq:hij},
we can construct a $\atba\smeover T
C_0(T)$-\eb{} as in \cite{ra1,90c}.  Specifically, we set
\begin{equation*}
\Y'=\set{(y_i)\in\prod_I\xittxib:\hij\(y_j^\fij\)=
y_i^\fij}.
\end{equation*}
{}From \eqref{eq:hij} we deduce that if $t\in\fij$
and $x=(x_i), y=(y_i)\in\Y'$, then
\begin{equation*}
\lip (\atba)^{F_i} <x_i,y_i>(t)=\lip (\atba)^{F_j}<x_j,y_j>(t).
\end{equation*}
(Since $\atba$ has Hausdorff spectrum $T$, we may view it as the
section algebra of a \cs-bundle over $T$.)
Since a similar equation holds for the $C(F_i)$-valued inner
products, we obtain
well-defined sesqui-linear forms on $\Y'$ by the formulas
\begin{equation}\label{eq:ips}
\begin{split}
\lip\atba<x,y>(t)&=\lip(\atba)^{F_i}<x_i,y_i>(t),\text{ and}\\
\rip C_0(T)<x,y>(t)&=\rip C(F_i)<x_i,y_i>(t) ,
\end{split}
\end{equation}
for $t\in F_i$.
Notice that $\Y'$ admits natural left and right actions of $\atba$ and
$C_0(T)$, respectively.  The next lemma can be proved along the
lines of \cite[Proposition~2.3]{ra1}.
\begin{lem}\label{lem:ymod}
With the inner products given by \eqref{eq:ips},
$$\Y=\set{y\in\Y':\text{$t\mapsto\rip C_0(T)<y,y>(t)$ vanishes at
infinity}}$$
is a complete $\atba\smeover T C_0(T)$-\eb.
\end{lem}

While $\Y$ is the sort of module required in \thmref{thm:crocker},
it unfortunately carries
no obvious $G$-action---let alone one equivalent to $\tau$.
To overcome this, we will want to see
that $\Y$ is isomorphic to a special subalgebra of $A$.  To do this
let
$$\n=\set{a\in A:\text{$t\mapsto\tr(a^*a)(t)$ is in $C_0(T)$}}.$$
Then $\rip C_0(T)<x,y>(t)=\tr(x^*y)(t)$ defines a $C_0(T)$-valued
inner product on $\n$ (\cite[4.5.2]{dix}).  Because $A$ has
continuous trace,
$\n$ is dense in $A$ by
Definition~4.5.2 and Lemma~4.5.1(ii) of \cite{dix}; thus
$\sp\set{\rip C_0(T)<x,y>:x,y\in\n}$
is an ideal in $C_0(T)$ without common zeros, and
hence is dense in $C_0(T)$.
The next result is a pleasant surprise.
\begin{lem}
With respect to the norm $\|a\|_2=\|\rip C_0(T)<a,a>\|_\infty^{1/2}$,
$\n$ is a (full) right Hilbert $C_0(T)$-module.
\end{lem}
\begin{pf}
The only issue is to see that $\n$ is complete.
Observe that if $a\in \n$, then $a(t)$ is a Hilbert-Schmidt
operator, and $\|\cdot\|_2$ induces the Hilbert-Schmidt norm
$\hs\|\cdot\|$ on
$\n(t)=\set{a(t)\in A(t):a\in \n}$.
In particular, for any $t\in T$,
$\hs\|a(t)\|\le\|a\|_2$.
So suppose that $\set{a_n}$
is $\|\cdot\|_2$-Cauchy in $\n$.  Since
the \cs-norm $\|\cdot\|$ is dominated by $\|\cdot \|_2$,
we have $a_n$ converging to some $a$ in $A$.  Since the
Hilbert-Schmidt operators on any Hilbert space are complete in
the Hilbert-Schmidt norm, the Cauchy sequence $a_n(t)$ must
converge, and must converge to $ a(t)$ in the Hilbert-Schmidt norm and
$\tr(a^*a)(t)<\infty$.

We still have to show that $a\in \n$ and that $a_n$ converges to
$a$ in $\n$.
Fix $\epsilon>0$.
Choose $N$ so that $n,m\ge N$ implies that $\|a_n-a_m\|_2<{
\epsilon / 2}$.
If $t\in T$, then there is a $k\ge N$ so that $\hs\|a_k(t)-a(t)\|<
{\epsilon / 2}$.  Then if $n\ge N$,
\begin{equation*}
\begin{split}
\hs\|a_n(t)-a(t)\|&\le\hs\|a_n(t)-a_k(t)\|+\hs\| a_k(t)-a(t)\|
\\
&\le \|a_n-a_k\|_2 + \hs\| a_k(t)-a(t)\| <\epsilon.
\end{split}
\end{equation*}
Our result follows as $\|x\|_2=\sup_t\hs\|x(t)\|$.
\end{pf}

We also need the following technical result.  It is a special
case of \cite[Lemma~1]{green4}.
\begin{lem}[Green]
\label{lem:green}
If $x,y\in\X_i$ and $t\in F_i$, then
$\tr\(\lip A^{F_i}<x,y>(t)\)=\rip C(F_i)<y,x>(t)$.
\end{lem}
\begin{pf}
This result is proved for $x=y$ in the second paragraph of the
proof of Theorem~2.15 in \cite{dpw3}.
Since everything in sight is trace-class, the general case
follows from the
usual polarization
%\footnote{OK, guys, we can scribe Phil's
%proof word for word from \cite{green4} if you prefer.}
identities: $4\lip A^{F_j} <x,y> =\sum_{k=0}^3\lip A^{F_j}
<x+i^ky, x+i ^k y>$ and $4\rip C_0(T)<x,y> = \sum_{k=0}^3 \rip
C_0(T) <i^k x + y, i^k x + y>$.
\end{pf}

We define a map $\Phi_i:\xittxib
\to\n^{F_i}$ as follows. Suppose
$y_i=\sum_k x_k\tensor \flat(z_k)$ is a sum of elementary
tensors in $\X_i\odot\bX_i$.
Then for $t\in F_i$,
$$\Phi_i(y_i)(t)=\sum_k\lip A^{F_i}<x_k,z_k>(t),$$
defines a map on
$\X_i\odot\bX_i$ (it is sesqui-linear),
which preserves inner products by the following
computation.  Let $y_i'=\sum_k u_k\tensor \flat(v_k)$.  Then
\begin{align}\label{eq:isometric}
\brip C(F_i)<\Phi_i(y_i),\Phi_i(y_i')>
&= \tr\(\Phi_i(y_i)^*\Phi(y_i')(t)\) \\
&= \tr \bigl(\sum_{k,r} \lip A^{F_i}<x_k,z_k>^* \lip A^{F_i}
<u_r,v_r> \bigr) (t) \notag\\
& = \tr \bigl( \sum_{r,k} \blip A^{F_i} < z_k \cdot {\rip C(F_i)
<x_k, u_r>}, v_r> \bigr) (t) \notag\\
\intertext{which, using \lemref{lem:green}, is}
&= \sum_{r,k} \brip C(F_i) <v_r,z_k\cdot{\rip C(F_i) <x_k,u_r>}>(t)
\notag\\
&= \rip C(F_i)<y_i,y_i'>(t)\notag
\end{align}
Thus,  $\Phi_i$
extends to a map on $\xittxib$ taking values in $\n^{F_i}$ since
the latter is complete.  Notice that we may replace $\X_i$ by
$\X_i^\fij$ in the above to obtain a similar map $\Phi_i^\fij$ into
$\n^\fij$, and that for $t\in \fij$ we have $\Phi_i(y)(t)
= \Phi_i^\fij(y^\fij)(t)$ for any $y\in\xittxib$.
Now suppose that $y=(y_i)\in\Y$, $t\in\fij$, and $\epsilon>0$.
Choose $\tilde y_j\in \X_j\odot\bX_j$ so that $\|y_j-\tilde y_j
\|<\epsilon$.  Thus $\|y_j^\fij -\tilde y_j^\fij\|<\epsilon$, and
\eqref{eq:isometric} implies that $|\Phi_j(y_j)(t)-\Phi_j( \tilde
y_j)(t)|<\epsilon$.  As $y\in\Y$, $\|y_i^\fij -\hij( \tilde
y_j^\fij)\|<\epsilon$, and a calculation on elementary tensors
shows that $\Phi_i
\(\hij(\tilde y_j)\)(t)=\Phi_j(\tilde y_j^\fij)(t)$.  It follows that
$|\Phi_i(y_i)(t)-\Phi_j(y_j)(t)|<\epsilon$; since $\epsilon$ was
arbitrary, we have $\Phi_i(y_i)(t)=\Phi_j(y_j)(t)$
Thus we can define $\Phi:\Y\to\n$ by setting $\Phi\((y_i)
\)(t)=\Phi_j(y_j)(t)$ for $t\in F_j$.  We have shown above that this
is well-defined, and it follows from \eqref{eq:isometric} that $\Phi$
does indeed take values in $\n$.
\begin{prop}\label{prop:phiiso}
The map $\Phi$ defined above extends to a Hilbert $C_0(T)$-module
isomorphism from $\Y$ {\em onto} $\n$.
\end{prop}
\begin{pf}
Since we have already shown that $\Phi$ preserves inner products,
and $\Phi$ is clearly $C_0(T)$-linear, we only have to show that
$\Phi(\Y)$ is dense in $\n$.
Observe that the $C_0(T)$-submodule $\Phi(\Y)$ is
also an ideal in $A$:
$a\Phi(y)=\Phi\((a\tensor 1)\cdot y\)$ and
$\Phi(y)a=\Phi\((1\tensor\flat(a^*))\cdot y\)$.
Since $\Phi(\Y)$ is certainly \cs-norm dense in $A$, we have that
$\Phi(\Y)(t)$ is norm dense in $A(t)$ for each $t\in T$.  In
particular, $\Phi(\Y)(t)$ contains the finite-rank operators.
(The finite-rank operators are the Pedersen ideal in $\K$
\cite[\S5.6]{ped}.)
Therefore $\Phi(\Y)(t)$ is dense in $X(t)$ in the Hilbert-Schmidt
norm.  Now fix $a\in\n$ and $\epsilon>0$.
Choose a compact set $C\subseteq T$ such that $\hs\|a\|^2=
|\tr(a^*a)(t)|<{\epsilon/ 4}$ if $t\notin C$.  For each $t\in C$,
there is a $y\in\Y$ such that $\hs\|\Phi(y)(t)-a(t)\|<{\epsilon/
4}$, and a relatively compact
\nbhd{} $U$ of $t$ such that
\begin{equation}
\label{eq:dagger}
\hs\|\Phi(y)(t')-a(t')\|<{\epsilon\over 2}
\text{ for all $t'\in U$.}
\end{equation}
Next choose a partition of unity $\set{\phi_i}_{i=1}^n\subseteq
C_c(T)$
subordinate to a
cover $U_1,\dots,U_n$ of $C$ as in \eqref{eq:dagger} for elements
$y_1,\dots,y_n$ in $\Y$.  (That is, $\sum_i\phi_i\equiv 1$ on $C$
and $\sum_i\phi_i\le1$ otherwise.)
Let $y=\sum y_i\cdot\phi_i$.  Since $\|a-a\cdot\sum\phi_i\|_2 <
{\epsilon/ 2}$, it follows that $\|\Phi(y)-a\|_2<\epsilon$.
This proves density and completes the proof.
\end{pf}
\begin{cor}
Let $\K(\n)$ denote the compact operators on the Hilbert
$C_0(T)$-module
$\n$.  (In less modern terms, $\K(\n)$ is the imprimitivity
algebra of the right $C_0(T)$-rigged space $\n$.)
Then there is a $C_0(T)$-isomorphism $Q$
of $\atba$ onto $\K(\n)$ which satisfies
\begin{equation}
Q\(a\tensor\flat(b)\)(c)=acb^*.
\end{equation}
\end{cor}
\begin{pf}
Recall from \lemref{lem:ymod} that $\atba\cong\K(\Y)$.
Define $Q:\L(\Y)\to\L(\n)$ by $Q(T)(x)=\Phi\(T\(
\Phi^{-1}(x)\)\)$.  Then check that
$Q\(\K(\Y)\)=\K(\n)$, and note that $Q$ is
$C_0(T)$-linear.

Finally, let $y=\(x_i\tensor\flat(y_i)\)$ be an element of $\Y$
whose components are elementary tensors, and
let $T$ be left-multiplication by $a\tensor\flat(b)$.  Then elements
of the form $\Phi(y)$ span a dense subset
of $\n$\footnote{To see this,
note that $\n(t)$ is a Hilbert space and a partition of unity
argument shows that it is enough to see that the span of such
elements is dense in each $\n(t)$.  Now the assertion follows
from the following: (1)~elementary tensors span a dense subset of
$\xittxib$; (2)~if $t\in F_i$, $\Phi(y)(t)=\Phi_i(y)(t)$; and
(3)~if $y_0\in\X_i\odot\bX_i$, then
the argument in \cite[Lemma~2.2]{ra1}
shows that there is a $y=(y_i)\in \Y$ with
$y_i=y_0$ and with the other $y_j\in\X_j\odot\bX_j$.}
and $\Phi\(T(y)\)=\Phi\((ax_i\tensor\flat(by_i))\)=a\Phi(y)b^*$.
\end{pf}

We now claim that, for each $s\in G$, $\alpha_s(\n)=\n$.
By assumption on $(A,\alpha)$, if $a\in A$ and
$\pi_t$ denotes evaluation at $t$, then
$\pi_{s^{-1}\cdot t}$ is equivalent to $\pi\circ\alpha_s$, and
\begin{equation*}
\alpha_s(a^*a)(t)
=\pi_t\(\alpha_s(a^*a)\) = V\pi_{s^{-1}\cdot t}(a^*a)V^*
= V a^*a(s^{-1}\cdot t)V^*
\end{equation*}
for some unitary $V$.
It follows that $\tr\(\alpha_s(a^*a)(t)\) = \tr \(a^*a(s^{-1}
\cdot t)\)$ for all $a\in A$.  Therefore if $x,y\in\n$, then a
polarization argument implies that
\begin{equation}
\label{eq:trace}
\tr\(\alpha_s(x^*y)(t)\)=\tr\((x^*y)(s^{-1}\cdot t)\),
\end{equation}
satisfying the claim.
\begin{prop}
The action defined by $u_s(x)=\alpha_s(x)$ is strongly continuous
on $\n$ and satisfies
$$\brip C_0(T)<u_r(x),u_r(y)>(t)=\rip C_0(T) <x,y>(r^{-1}\cdot
t)$$
for all $x,y\in \n$, $t\in T$, and $s\in G$.
\end{prop}
\begin{pf}
The second assertion follows from \eqref{eq:trace}.  Thus,
because we
have already shown that $\alpha_s(\n)=\n$, we only have to show strong
continuity.

We first claim that $\m=\n^2$ is $\|\cdot\|_2$-norm dense in $\n$.
(In the notation of \cite{dix}, $\n$ coincides with $\frak n$.)
By \cite[4.5.1]{dix}, $\m$ is \cs-norm dense in $A$, so $\m(t)$
contains the finite-rank operators for each $t\in T$ and hence is
dense in $\n(t)$ in the $\hs\|\cdot\|$-norm.
Thus given $t_0\in T$ and
$\epsilon>0$, there is a \nbhd{} $U$ of $t_0$ and  $b
\in \m$ such that $\hs\|a(t)-b(t)\|<{\epsilon/ 2}$ for all
$t\in U$.  Another partition of unity argument
as in \propref{prop:phiiso} implies that there is a $b\in\m$ such
that $\|a-b\|_2<{\epsilon }$; this establishes the claim.

Since each $u_t$ is $\|\cdot\|_2$-isometric, it suffices to show that
$\lim _{s\to e}\|\alpha_s(a)-a\|_2=0$; from the previous paragraph we
can assume that $a\in\m$.
It follows from \eqref{eq:trace} that $a(t)$ and $\alpha_s(a)(t)$
are trace class operators.
Recall that $\|T\|_1=\tr\(|T|\)$ is a norm on the trace-class
operators.  While it is apparent that
$\|a(t)\|_1=\|\alpha_s(a)(s\cdot t)\|_1$ if $a\ge0$ (see
\eqref{eq:trace}), we do not know whether this holds in general as it
is not obvious that $|a|\in\m$ if $a$ is. However, $a=\sum_{i=1}^n
\alpha_ia_i$ where $\alpha_i\in\C$ and $a_i\in\m^+$
by
\cite[4.5.1]{dix}.
Thus,
$\|\alpha_s(a)(t)\|_1\le\sum_i|\alpha_i|\|\alpha_s(a_i)\|_1=\sum_i
|\alpha_i|\|a_i(s^{-1}\cdot t)\|_1$, and there is a constant $K$,
which depends only on $a$ (and not on  $s\in G$),
so that $\sup_t \|\alpha_s(a)(t)\|_1\le K$.
Therefore, if $\|\cdot \|$ denotes the \cs-norm, then
\begin{align*}
\|\alpha_s(a)-a\|_2 &=\sup_t\|(\alpha_s(a)-a)^*(\alpha_s(a)-a)(t)\|_1\\
&\le \sup_t\|\alpha_s(a)-a\|\,\|(\alpha_s(a)-a)(t)\|_1 \\
&\le \|\alpha_s(a)-a\|\(K+\sup_t\|a(t)\|_1\),
\end{align*}
where the second inequality follows from \cite[3.4.10]{pednow}.
The conclusion now follows from the strong continuity of $\alpha$
in the \cs-norm on $A$.
\end{pf}

\begin{lem}[{\cite[\S3]{combes}}]
Suppose that $\X$ is an $A\sme B$-\eb{}
and that $u:G\to\Aut(\X)$
is an action of $G$ on $\X$: that is, there is a strongly
continuous automorphism $\tau:G\to \Aut(B)$ such that $\tau_s\(
\rip B<x,y>\cdot b\)=\brip B<u_s(x),u_s(y)>\cdot\tau_s(b)$.
Then, if $T\in\L(\X)$,
\begin{equation}\label{eq:above}
\alpha_s(T)=u_sTu_s^{-1}
\end{equation}
is in $\L(\X)$, and \eqref{eq:above} defines
a strongly continuous automorphism group
$\alpha:G\to\Aut(A)$ satisfying $\alpha_s\(a\cdot\lip A<x,y>\)=
\alpha_s(a)\cdot\blip A<u_s(x),u_s(y)>$.
\end{lem}
\begin{pf}
Certainly, $\alpha_s$ defines an automorphism of $\L(\X)$, and
since $u_s\lip<x,y>u_s^{-1}=\blip A<u_s(x),u_s(y)>$, $\alpha_s$
restricts to an automorphism of $A$.  The rest is
straightforward.
\end{pf}

\begin{pf*}{Proof of \thmref{thm:crocker}}
We have shown above that $(\n,u)$ is a Morita equivalence between
the systems $(C_0(T),\tau)$ and $(\atba,\beta)$ where $\beta_s
(a\tensor\flat(b))=u_s\(Q\(a\tensor\flat(b)\)\)u_s^{-1}$.
Now let $c\in\n$ and recall that $u_t$ is the restriction of $\alpha_t$
to $\n$.  Then
\begin{align*}
\beta_t\(a\tensor\flat(b)\)\cdot c
&=\alpha_t\(Q\(a\tensor\flat(b)\)\(\alpha_t^{-1}(c)\)\)
= \alpha_t(a)c\alpha_t(b)^* \\
&=Q\(\alpha_t (a)\tensor \bar\alpha_t(\flat(b))\)(c).
\end{align*}
Thus, $\beta=\alpha\tensor\bar\alpha$, and we have shown that
$[A,\alpha]^{-1}$ exists and equals $[\bA,\bar \alpha]$.  This
completes the proof.
\end{pf*}
%
%
%   SCCS Info: ckrw5.tex version 3.1: 10/19/93
%
\section{The Forgetful Homomorphism}
\label{sec:5}

In this section we will require that $H^2(T;\Z)$ be countable,
and as before, that $(G,T)$ be a second countable locally compact
transformation group.  The homomorphism $F:\brgt\to\brt$ defined
by $F\([A,\alpha]\)=\delta(A)$ is called the {\it Forgetful
Homomorphism} (where we identify $\brt$ with $H^3(T;\Z)$ as in
Remark~\ref{rem:3.3}). The image of $F$ is of considerable interest
as it describes exactly which stable algebras admit actions inducing
a given action on $T$.  More precisely, we have the following lemma.
\begin{lem}
If $A$ is a
stable, separable continuous-trace \cs-algebra with spectrum $T$, then
$A$ admits an an automorphism group $\alpha:G \to\Aut(A)$ inducing the
given action on $T$ if and only if $\delta(A)$ is in $\Im(F)$.
\end{lem}
\begin{pf}
{}From the definitions it is clear that $\delta(A)$
is in the image of $F$ if and only if $A$ is (strongly) Morita
equivalent over $T$ to an algebra $B$ which admits an action
$\alpha$ covering $\tau$; i.e., $(B,\alpha)\in\sbrgt$.  But then
$(B\tensor\K,\alpha\tensor\id)\in\sbrgt$, and $A\tensor\K\cong
B\tensor\K$ by \cite{bgr}, and we have to check that this
isomorphism is $C_0(T)$-linear.
But since $\X$ is an $A\smeover T B$-bimodule, there is a natural
action of $C_0(T)$ on the linking algebra $C$; since the
isomorphism of $B\tensor\K$ and $A\tensor\K$ with $C\tensor \K$ are
obtained by conjugation by partial isometries in $\M(C\tensor
\K)$ \cite[\S2]{bgr}, they are $C_0(T)$-linear.
%The lemma will follow once we
%verify that all the stable isomorphisms involved can be taken
%to be $C_0(T)$-linear.
%First note that the isomorphism
%$A\tensor\K\cong
%B\tensor\K$ induces a homeomorphism of $T$ which takes
%$\delta(A)$ to $\delta(B)$.  Moreover,
%$\delta(A)=\delta(B)$ since these algebras are Morita equivalent
%over $T$ (see, for example, \cite[Theorem~3.5]{90c}).  Thus by
%\cite[Theorem~2.22]{pr1}, there is an automorphism of $B\tensor
%\K$ inducing the inverse homeomorphism on $T$.  Composition gives
%an isomorphism of $A\tensor\K$ with $B\tensor \K$ which induces
%the identity on $T$, and is therefore $C_0(T)$-linear by
%\cite[Theorem~2.22]{pr1} again.
Finally, if $A$ is stable, then
$A$ is $C_0(T)$-isomorphic to $A\tensor\K$ by \cite[Lemma~4.3]{pr2}.
This proves the lemma.
\end{pf}

As an example of the significance of these ideas, notice that
\cite[Theorem~4.12]{rr} implies that $F$ is surjective when $G=\R$
and $(G,T)$ is a principal $\T$-bundle (provided, say, $T/G$ is a
$CW$-complex).  The analysis of this section will give a
substantial generalization of this result.
Our approach is to identify three obstructions to an element $\c\in
H^3(T;\Z)$ being in $\Im(F)$, and then to show that the vanishing
of these obstructions is sufficient (as well as necessary).

The first obstruction is that in order that $\c\in H^3(T;\Z)$
be in $\Im(F)$, we must have
\begin{equation} \label{eq:5.1}
\c\in H^3(T;\Z)^G=\set{\c\in H^3(T;\Z): \text{$s\cdot\c=\c$ for
all $s\in G$}}
\end{equation}
(Here and in the sequel, we view $H^n(T;\Z)$ as a $G$-module via
the $G$-action on $T$; that is, if $\ell_s$ denotes the
homeomorphism $t\mapsto s\cdot t$ of $T$, then $s\cdot\c=\ell_{
s^{-1}}^*(\c)$.)
The necessity of \eqref{eq:5.1}
is a consequence of the following lemma which, although we
present a different proof here, is contained in Theorem~2.22 of
\cite{pr1}.
\begin{lem}
Suppose that $A$ is a continuous-trace \cs-algebra with spectrum
$T$, that $\alpha\in\Aut(A)$, and that $h$ is the homeomorphism
of $T$ induced by $\alpha$.  Then $h^*\(\delta(A)\)=\delta(A)$.
\end{lem}
\begin{pf}
Let $h^*(A)=C_0(T)\tensor_{C_0(T)}A$ be the pull-back of $A$
along $h:T\to T$.  Then $\delta\(h^*(A)\)=h^*\(\delta(A)\)$ by
\cite[Proposition~1.4(1)]{rw}.  On the other hand $\phi\tensor
a\mapsto\phi\tensor\alpha(a)$ extends to a $C_0(T)$-isomorphism
of $\id^*(A)$ onto $h^*(A)$.  The result follows.
\end{pf}

Our other obstructions are obtained by identifying a subgroup
of the Moore group $H^3\(G,C(T,\T)\)$ and defining homomorphisms
$d_2$ of $H^3(T;\Z)^G$ into $ H^2\(G,H^2(T;\Z)\)$ and $d_3$ from
the kernel of $d_2$ to the corresponding quotient of
$H^3 \(G,C(T,\T)\)$.  We then show that $\c\in\Im(F)$ if and only
if $\c\in H^3(T;\Z)^G$ and $d_2(\c)=d_3(\c)=0$.
This will occupy the remainder of this section.

The basic idea for the construction of $d_2$ is as follows.
Let $\ell:G\to\homeot$ be the map induced by $(G,T)$.  Notice
that if $\c\in H^3(T;\Z)^G$, then $\ell(G)\subseteq \Homeo_\c(T)$.
If $A$ is the essentially unique {\it stable\/} continuous-trace
\cs-algebra with $\delta(A)=\c$, then by
\cite[Theorem~2.22]{pr1} there is a short exact sequence
{\dgARROWLENGTH=1em
\begin{ediagram*}
\node{1} \arrow{e}
\node{\Autc (A)/\Inn(A)} \arrow{e}
\node{\Out(A)} \arrow{e} \node{\Homeo_\c(T)}
\arrow{e}\node{1.}
\end{ediagram*}}%
Furthermore, it follows from \cite[Theorem~2.1]{pr1} that there
is an isomorphism $\zeta:\Autc (A)/\Inn(A)\to H^2(T;\Z)$.
Therefore there should be an obstruction $d_2(\c)$ in $H^2 \(G,
H^2(T;\Z)\)$ to lifting $\ell$ to a homomorphism $\gamma: G \to
\Out(A)$.

The existence of $d_2$ will follow from the next lemma.  Notice
that if $N$ is a normal {\it abelian\/} subgroup of a group $H$,
then $H/N$ acts on $N$ by conjugation.
\begin{lem}\label{lem:first}
Suppose that $H$ is a Polish group, that $N$ is a closed normal
abelian subgroup, and that $\ell:G\to H/N$ is a continuous
homomorphism.
Then $N$ is a $G$-module (where $g\in G$ acts on $n\in N$ by $g
\cdot n=\ell_gn\ell_g^{-1}$), and there is a cohomology class
$\c(\ell)\in H^2(G,N)$ which vanishes if and only if there is a
continuous homomorphism $\gamma:G\to H$ which lifts $\ell$ (i.e.,
$\gamma_gN=\ell_g$).  In fact, one obtains a cocycle $n\in
Z^2(G,N)$ representing $\c$ by taking any Borel lift $\gamma'$
of $\ell$, and defining $n$ by
\begin{equation}
\gamma_s'\gamma_t'=n(s,t)\gamma_{st}'.
\label{eq:1}
\end{equation}
\end{lem}
\begin{pf}
By \cite[Proposition~4]{moore3} we can find a Borel section $s:H/N
\to H$ for the quotient map such that $s(N)=e$.  Define
$\gamma'=s\circ\ell$, and define $n$ by \eqref{eq:1}.
Then $n$ is certainly Borel and comparing
$\gamma_r'(\gamma_s'\gamma_t')$ with $(\gamma_r'\gamma_s')\gamma_t'$
shows that  $n$ is
indeed a cocycle.
A standard argument shows that the class of $n$ in $H^2(G,N)$ is
independent of our choice of section $s$.

Evidently, if $\gamma'$ is a homomorphism, then $[n]=0$ as $n$ is
identically equal to $e$.  On the other hand, if $[n]=0$,
then there is a Borel function $\l:G\to N$ such that $\l_e=e$ and
$n(s,t)=(s\cdot\l_t)^{-1}\l_s^{-1}\l_{st}$,
and
\begin{equation*}
(\l_s\g_s')(\l_t\g_t') = \l_s(s\cdot\l_t)\g_s'\g_t'
= \l_s(s\cdot\l_t) n(s,t)\g_{st}'
= \l_{st}\g_{st}'.
\end{equation*}
Therefore $\g=\l\g'$ is a Borel homomorphism lifting $\ell$,
which is continuous by \cite[Proposition~5]{moore3}.
\end{pf}

In order to apply \lemref{lem:first} we have to see that the
groups involved are Polish.  However, because $A$ is a separable
\cs-algebra, then $\Aut(A)$, with the topology of pointwise
convergence, is a Polish group.  Then, as $\Autc (A)$ is
closed in $\Aut(A)$, it is also Polish.  Since we are assuming
that $H^2(T;\Z)$ is countable, $\Inn(A)$ is open in
$\Autc (A)$ and closed in $\Aut(A)$ by
\cite[Theorem~0.8]{rr}.  In particular,
$\Out(A)$ is Polish, as is $H^2(T;\Z)\cong
\Autc (A)/\Inn(A)$ which is even discrete.
Finally we give $\homeot$ the compact open topology.  Then it is
not hard to see that the map $\rho:\Aut(A)\to\homeot$ is
continuous.  In fact, $\homeot$ is homeomorphic to $\Aut\(C_0(T)\)$,
and so is certainly a Polish group as well.  We have not
been able to show that $\Homeo_{\delta(A)}(T)$ is closed in $\homeot$, so
we don't know for sure that it is Polish.  But,
as it follows from \cite[Theorem~2.22]{pr1} that
$\rho$ defines a continuous
injection $h:\Aut(A)/\Autc (A)\to\homeot$ with image
exactly $\Homeo_{\delta(A)}(T)$,  Souslin's Theorem
\cite[Theorem~3.3.2]{arveson} implies that
$\Homeo_{\delta(A)}(T)$ is a Borel subset of $\homeot$ and $h$ is
a Borel isomorphism.
Thus we can view $\ell$ as a Borel homomorphism of $G$ into the
Polish group $\Out(A)/H^2(T;\Z)$, which is automatically continuous.
We can now apply \lemref{lem:first} to
{\dgARROWLENGTH=1.5em
\begin{ediagram*}
\node{0}\arrow{e}\node{H^2(T;\Z)}\arrow{e}
\node{\Out(A)}\arrow{e}\node{\Out(A)/H^2(T;\Z)}\arrow{e}
\node{0} \\
\node[4]{G}\arrow{nw,b,..}{\alpha} \arrow{n,r}{\ell}
\end{ediagram*}}%
to get the desired obstruction $d_2\(\delta(A)\)$ in
$H^2\(G,H^2(T;\Z)\)$.  However there is a small point to check.
The $G$-action on $H^2(T;\Z)$ coming from \lemref{lem:first} is
that coming from the identification of $H^2(T;\Z)$ with the
subgroup $\Autc(A)/\Inn(A)$ of $\Out(A)$
via the isomorphism $\zeta$.  The $G$-action,
then, is induced by conjugation by elements of $\Aut(A)/\Autc{A}
\cong_h \Homeo_{\delta(A)}(T)$.  Thus the action
of $s\in G$ on $[\phi]\in\Aut(A)/\Autc{A}$ is given by
$s\cdot[\phi]=[\gamma_s\phi\g_s^{-1}]$ where $\rho(\g_s)=\ell(s)$.
On the other hand, the usual action of $s\in G$ on $\c\in
H^2(T;\Z)$ is given by $s\cdot \c=\ell_{s^{-1}}^*(\c)$.
It is a relief that these actions
coincide.
\begin{lem}\label{lem:5.3}
Let $\zeta:\Autc(A)/\Inn(A)\to H^2(T;\Z)$ be the isomorphism from
\cite[Theorem2.1]{pr1}.  If $\g\in\Aut(A)$ and $\phi\in\Autc(A)$,
then
\begin{equation*}
\zeta\([\g\phi\g^{-1}]\)=h(\g^{-1})^*\(\zeta[\phi]\).
\end{equation*}
\end{lem}
\begin{pf}
To define $\zeta$, we follow \cite[\S5]{russell}.  View $A$ as
the sections of a \cs-bundle $\xi$ over $T$.
Then $\phi$ is locally implemented by multipliers.  Thus there are
an open cover $\set{N_i}$ of $T$ and $u_i\in\M(A)$ such that
$\phi(a)(x)=u_i(x)a(x)u_i^*(x)$ for $x\in N_i$.  (Recall that $u_i$
can be viewed as a field of multipliers in $\M(A(x))$ \cite{lee1}.)
Then $\zeta(\phi)$ is represented by the cocycle $\set{\l_{ij}}$
where $\l_{ij}(x)u_j(x)=u_i(x)$ for $x\in N_{ij}$.

Let $u_i$ and $\l_{ij} $ be as above.  For notational
convenience, let $y= h(\g)^{-1}(x)$.
Define an isomorphism $\g_x$
from $A(y)\to A(x)$ by $\g_x\(a(y)\)
=\g(a)(x)$.  Then for $x\in h(\g)(N_i)$, we have
\begin{align*}
(\g\phi\g^{-1})(a)(x)
&= \g_x\(\phi\g^{-1}(a)(y)\) \\
&= \g_x\(u_i(y)\g^{-1}(a)(y)
u_i(y)^*\) \\
&= \Ad\(\g_x\(u_i(y)\)\)\bigl[\g_x\(\g^{-1}(a)(
y)\)\bigr] \\
&= \Ad\(\g_x\(u_i(y)\)\)\bigl[a(x)\bigr]=\Ad\(\g(u_i)(x)\)[a(x)],
\end{align*}
so $\g\phi\g^{-1}$ is implemented over $h(\g)(N_i)$ by
$v_i=\g(u_i)$.  Therefore for $x\in h(\gamma)(N_{ij})$, we have
\begin{equation*}
\l_{ij}(y)v_j(x) =\l_{ij}(y)\g(u_j)(x)
= \g_x(\l_{ij}(y)u_j(y)\) =\g_x\(v_i(y)\)=
v_i(x).
\end{equation*}
Thus $\zeta(\g\phi\g^{-1})$ is represented by the cocycle
$\set{h(\g)(N_i),\l_{ij}\circ h(\g)^{-1}}$, which also represents
$h(\g^{-1})^* \(\zeta(\phi)\)$.  This completes the proof.
\end{pf}
\begin{lem}
The map $d_2$ defined above is a homomorphism from $H^3(T;\Z)^G$ to
$H^2\(G,H^2(T;\Z)\)$.
\end{lem}
\begin{pf}
Let $A$ and $B$ be stable continuous-trace \cs-algebras with
spectrum $T$.  Since there is a Borel section for the quotient
map from $\Aut(A)$ onto $\Out(A)$, we may assume that there are
Borel maps $\gamma:G\to\Aut(A)$ and $\delta:G\to\Aut(B)$ so that
the obstructions are determined, respectively,
by Borel maps
$n:G\times G\to\Autc (A)$ and $m:G\times G\to \Autc B$
such that $\gamma_s\g_t=n(s,t)\g_{st}$ and $\delta_s\delta_t =
m(s,t)\delta_{st}$.  (Here, $m$ and $n$ need not be
cocycles---only their images in $\Autc (A)/\Inn(A)$ and
$\Autc (B)/\Inn(B)$.)
Since $\g_s$ and $\delta_s$ induce the same homeomorphism of $T$,
$\g_s\tensor \delta_s$ defines an automorphism of $C=A\Ttensor B$
which also induces the same homeomorphism of $T$.
Moreover,
\begin{equation*}
(\g_s\tensor\delta_s)(\g_t\tensor\delta_t)=\( n(s,t)\tensor
m(s,t)\)(\g_{st}\tensor\delta_{st}),
\end{equation*}
so $d_2\(\delta(C)\)$ is determined by the cocycle
$[n\tensor m]$ in $H^2\(G,\Autc(C)/\Inn(C)\)$.
But under the isomorphism $\zeta_C:\Outc{C}\to H^2(T;\Z)$ we have
\begin{equation*}
\zeta_{C}(n\tensor m)=\zeta_A(n)+\zeta_B(m)
\end{equation*}
by \cite[Proposition~3.10]{pr2}.  Therefore, using \propref{prop:2.2},
we have
$d_2\(\delta(A)+\delta(B)\)=d_2\(\delta(A)\)+d_2\(\delta(B)\)$ as
required.
\end{pf}

We now turn to the definition of $d_3$.  The main technical tool
will be the following lemma.
Here we will need the twisted actions of
\cite{bs,para1,para3}: a twisted action of $G$ on
$A$ is a pair $(\alpha,u)$ consisting of Borel maps
$\alpha:G\to\Aut(A)$ and $u:G\times G\to\UM(A)$ satisfying the
axioms of \cite[Definition~2.1]{para1}.

\begin{lem}\label{lem:5.5}
Suppose that $A$ and $B$ are separable continuous-trace
\cs-al\-ge\-bras with spectrum $T$.
\begin{enumerate}
\item
Suppose that $\g:G\to\Out(A)$ is a continuous homomorphism and
that $C(T,\T)$ has the $G$-module
structure
coming from the
$G$-action induced  by $\g$ and the natural map
$h:\Out(A)\to\homeot$.
Then there is a class $d_A(\g)$ in $H^3\(G,C(T,\T)\)$ which
vanishes if and only if there is a twisted action, $\alpha:G\to
\Aut(A)$ and $u:G\times G\to\UM(A)$, such that
\begin{ediagram*}
\node{\Aut(A)}\arrow{e,t}{q} \node{\Out(A)} \\
\node[2]{G}\arrow{nw,b,..}{\alpha}\arrow{n,r}{\gamma}
\end{ediagram*}
commutes.
\item
Suppose that $\epsilon:G\to\Out(B)$ satisfies
$h_B\circ\epsilon=h_A\circ \g$.  Then there is
a continuous homomorphism $\g\tensor\epsilon:G\to \Out(A \Ttensor
B)$, and $d_{A\Ttensor
B}(\g\tensor\epsilon)=d_A(\g)d_B(\epsilon)$.
\item
The map $\alpha\mapsto\bar\alpha$ from $\Aut(A)\to\Aut(\bA)$, defined
by $\bar\alpha\(\flat(a)\)=\flat\(\alpha(a)\)$, maps $\Ad u$ to
$\Ad\flat(u)$, and hence induces an isomorphism of $\Out(A)$ onto
$\Out(\bA)$.  If  $\g:G\to\Out(A)$ is a continuous
homomorphism, and  $\bar\g:G\to\Out(\bA)$ is the corresponding map,
then
$d_{\bA}(\bar\g)=d_A (\g)^{-1}$.
\item If $\phi:A\to B$ is a
$C_0(T)$-isomorphism, then conjugation by $\phi$ induces a homomorphism
$\Ad(\phi):\Out(A)\to\Out(B)$ and $d_B\(\Ad(\phi)\circ\g\)=d_A(\g)$.
\end{enumerate}
\end{lem}
\begin{pf}
Choose a Borel lifting $\alpha$ of $\g$ such that $\alpha_e=\id$.
(Recall that we are assuming that $H^2(T;\Z)$ is countable so that
$\Out(A)$ is Polish.)  Since $q(\alpha_s\alpha_t)=q(\alpha_{st})$,
there is a Borel map $i:G\times G\to\Inn(A)$ such that
$\alpha_s\alpha_t= i(s,t)\alpha_{st}$, and $i(s,e)=i(e,t)=\id$ for all
$s,t\in G$. Since $\UM(A)\to\Inn(A)$ is a surjective homomorphism of
Polish groups, there is a Borel map $u:G\times G\to \UM(A)$ such that
\begin{equation}\label{eq:dagone}
\alpha_s\alpha_t=\Ad\(u(s,t)\)\circ\alpha_{st},\text{ and }
u(s,e)=u(e,t)=1. \end{equation}
Since $(\alpha_r\alpha_s)\alpha_t=\alpha_r(\alpha_s\alpha_t)$, we
must have
\begin{equation}
\Ad\(u(r,s)u(rs,t)\)=\Ad\(\alpha_r\(u(s,t)\)u(r,st)\).
\end{equation}
Therefore there is a Borel function $\l:G\times G\times G\to
\UZM(A)=C(T,\T)$ such that
\begin{equation}
\l(r,s,t)u(r,s)u(rs,t)=\alpha_r\(u(s,t)\)u(r,st).
\label{eq:deflam}
\end{equation}
Clearly, $\l(e,s,t)=\l(s,e,t)=\l(s,t,e)=1$ for all $s,t\in G$.
We want to show that $\l$ is a 3-cocycle for the action of $G$ on
$C(T,\T)$ defined above.  Notice that, because $\alpha_s$
lifts $\g_s$, if we view $C(T,\T)$ as the center of $\UM(A)$,
then the action of $G$ on $C(T,\T)$ is given
by $s\cdot f=\alpha_s(f)$.
Our computations follow \cite[\S IV.8]{maclane}.  Let $L=
\alpha_p\(\alpha_r\(u(s,t)\)u(r,st)\)u(p,rst)$.
Then on the one hand
\begin{align*}
L&=\alpha_p\(\l(r,s,t)u(r,s)u(r,st)\)u(p,rst) \\
&=\alpha_p\(\l(r,s,t)\)
\bigl[\l(p,r,s)u(p,r)u(pr,s)u(p,rs)^*\bigr] \\
&\qquad\bigl[\l(p,rs,t)u(p,rs)u(prs,t)u(p,rst)^*\bigr]
u(p,rst) \\
&=
\alpha_p\(\l(r,s,t)\)\l(p,r,s)\l(p,rs,t)u(p,r)u(pr,s)u(prs,t),\\
\intertext{while on the other hand,}
L&=
u(p,r)\alpha_{pr}\(u(s,t)\)u(p,r)^*\alpha_p\(u(r,st)\)u(p,rst)\\
&=
u(p,r)
\bigl[\l(pr,s,t)u(pr,s)u(prs,t)u(pr,st)^*\bigr]
u(p,r)^* \\
&\qquad\bigl[\l(p,r,st)u(p,r)u(pr,st)u(p,rst)^*\bigr]
u(p,rst) \\
&=\l(pr,s,t)\l(p,r,st)u(p,r)u(pr,s)u(prs,t).
\end{align*}
It follows that $\l$ is a 3-cocycle.

Next we observe that the class of $\l$ depends only on $\gamma$
and not on our choice of $u$ or $\alpha$.  First suppose that $v$
is another unitary-valued Borel map on $G\times G$ such that $\Ad
u=\Ad v$.  Then there is a Borel function $w:G\times G\to
C(T,\T)$ such that $v(s,t)=w(s,t)u(s,t)$.  Let $\mu$ be the
3-cocycle corresponding to $v$ as in \eqref{eq:deflam}.
Then
\begin{equation*}
\mu(r,s,t)w(r,s)w(rs,t)=\alpha_r\(w(s,t)\)w(r,st)\l(r,s,t),
\end{equation*}
so $\l$ and $\mu$ define the same class in $H^3\(G,C(T,\T)\)$.
If $\beta$ is another lift of $\g$, then there is a unitary
valued Borel function $\tv$ on $ G$ such that $\beta_s =
\Ad(\tv_s)\circ\alpha_s$.  Then we can choose the lift
\begin{equation*}
v(s,t)=\tv_s\alpha_s(\tv_t)u(s,t)\tv_{st}^*,
\end{equation*}
so that $\beta_s\beta_t=\Ad\(v(s,t)\)\beta_{st}$,
and compute that
\begin{align*}
\beta_r\(v(s,t)\)v(r,st)
&=\tv_r\alpha_r\(\tv_s\alpha_s(\tv_t)u(s,t)\tv_{st}^*\)\tv_r^*\tv_r
\alpha_r(\tv_{st})u(r,st)\tv_{rst}^* \\
&=\tv_r\alpha_r(\tv_s)\alpha_r\(\alpha_s(\tv_t)\)\alpha_r(u(s,t)
\) u(r,st)\tv_{rst}^*\\
&=\tv_r\alpha_r(\tv_s)
\bigl[u(r,s)\alpha_{rs}(\tv_t)u(r,s)^*\bigr]
\bigl[\l(r,s,t)u(r,s)u(rs,t)\bigr]
\tv_{rst}^* \\
&=\l(r,s,t)\tv_r\alpha_r(\tv_s)u(r,s)\tv_{rs}^*\tv_{rs}
\alpha_{rs}(\tv_t)u(rs,t)\tv_{rst}^* \\
&= \l(r,s,t)v(r,s)v(rs,t).
\end{align*}
Thus we get the same cocycle $\l$ for $\beta$ provided we choose
$v$ as above;  but since we have already observed that the class
of $\l$ in independent of our choice of $v$, we can conclude that
the class $d_A(\g)$ depends only on $\g$, as claimed.

If $d_A(\g)=0$, then $\l=\partial\mu$, and we can replace $u$ by $\mu
u$. (Then, of course, $\Ad u$ is unchanged and the corresponding $\l$
is identically one.)  Then it follows
from \eqref{eq:dagone} and \eqref{eq:deflam} that
$(\alpha,u)$ is a twisted action.
On the other hand, if there is a twisted action $(\alpha,u)$,
then the cocycle is certainly trivial.  This proves~(1).

Let $\g$, $\alpha$, $u$, and $\l$ be as above, and choose $\beta$
lifting $\epsilon$ as well as $v$ and $\mu$ in analogy with
\eqref{eq:dagone} and \eqref{eq:deflam}.  Since
$h_A\circ\g=h_B\circ \epsilon$, $\alpha\tensor
\beta$ defines a Borel map into $\Aut(A\Ttensor B)$.
Moreover,
\begin{align}
(\alpha_s\tensor\beta_s)(\alpha_t\tensor\beta_t)
&= \Ad\(u(s,t)\tensor v(s,t)\)\circ(\alpha_{st} \tensor \beta_{st}
) \label{eq:3} \\
\intertext{and with $w(s,t)=u(s,t)\tensor v(s,t)$, we have}
(\alpha_r\tensor\beta_r)\(w(s,t)\)w(r,st) &= \l\mu(r,s,t) w(r,s)
w(rs,t).\label{eq:4}
\end{align}
Then \eqref{eq:3} implies that $\alpha\tensor \beta$ defines a
Borel, hence continuous, homomorphism $\g\tensor\epsilon$ of $G$
into $\Out(A\Ttensor B)$,
satisfying $h_{A\Ttensor B}\circ \g\tensor
\epsilon=h_A\circ\g=h_B\circ\epsilon$,
and \eqref{eq:3} and \eqref{eq:4} together imply that $d_{A\Ttensor
B}(\g\tensor\epsilon)=d_A(\g)d_B(\epsilon)$.  This proves~(2).

Part (3) is easy: $\M(\bA)$ is naturally isomorphic to
$\overline{\M(A)}$, $\bar\alpha$ is a lift of $\bar\g$, and $\flat(u)$
satisfies $\bar\alpha_s\bar\alpha_t =
\Ad\(\flat\(u(s,t)\)\)\bar\alpha_{st}$.  But by definition of $\bA$,
applying $\flat$ to \eqref{eq:deflam} replaces $\l$ by
$\bar\l$, which is the inverse of $\l$ in $H^3$.

Finally, if $\alpha$ is a lift of $\g$, then $\beta_s=\phi\circ
\alpha_s\circ\phi^{-1}$ is a lift of $\Ad(\phi)\circ\g$.  But
then $\beta_s\beta_t=\Ad\(\phi\(u(s,t)\)\)[\beta_{st}]$.  Since
$\phi$ is $C_0(T)$-linear, the obstruction to $\phi\(u(\cdot,\cdot
)\)$ being a cocycle is that same as that for $u$.  That completes
the proof of the lemma.
\end{pf}

To define $d_2':H^1\(G,H^2(T;\Z)\)\to H^3\(G,C(T,\T)\)$ we
apply the above lemma to $(A,\alpha)=\(C_0(T,\K),\tau)$.
 Recall that $\zeta=\zeta_A:\Autc(A)/\Inn(A)\to H^2(T;\Z
)$ is a $G$-equivariant isomorphism (\lemref{lem:5.3}).
Thus, if $\rho\in Z^1\(G,H^2(T;\Z)\)$, then
$\rho'_s=\zeta^{-1}(\rho_s)\circ\ttau_s$ defines a continuous map
(\cite[Theorem~3]{moore3}) of $G$ into $\Out(A)$ such that $h\circ
\rho'=\ell$.  (We will abuse notation slightly and use
$\zeta^{-1}(\rho)$ to denote a representative in $\Autc(A)$ of the
class  $\zeta^{-1}(\rho)$ in $\Autc(A)/\Inn(A)$.)  Thus,
\lemref{lem:5.5}(1) gives us a class $d_A(\rho')$ in
$H^3\(G,C(T,\T)\)$ with the action of $G$ on $C(T,\T)$ being the
expected one: $s\cdot f(t) = f\(\ell_{s^{-1}}(t)\)=f(s^{-1}\cdot t)$.
If $\rho,\sigma\in Z^1\(G,H^2(T;\Z)\)$, then as $h\circ\rho'=h\circ
\sigma'$, part~(2) of our lemma implies that $d_{A\Ttensor A}(\rho'
\tensor\sigma')=d_A(\rho')d_A(\sigma')$.  But there is a
$C_0(T)$-isomorphism $\phi$ of $A\Ttensor A$ onto $A$.  Therefore
\begin{align*}
d_{A\Ttensor A}(\rho'\tensor\sigma')
&= d_A\(\phi\circ
(\rho'\tensor\sigma')\circ\phi^{-1}\) \\
&=
d_A(\phi\circ\zeta_{A\Ttensor A}^{-1}( \rho\sigma) \circ \phi^{-1})
\\
&= d_A(\zeta_A^{-1}\circ \rho\sigma)=d_A\((\rho\sigma)'\);
\end{align*}
the second equality is a consequence of
\cite[Proposition~3.10]{pr2}, the third follows from the general
fact that if $\phi:A\to B$ is a $C_0(T)$-isomorphism of
continuous-trace \cs-algebras, then
$\zeta_B\(\Ad(\phi)[\theta]\)=\zeta_A(\theta)$, which results from the
observation that if $\theta $ is locally implemented by $w_i\in
\M(A)$, then $\Ad(\phi)[\theta]$ is implemented by $\phi(w_i)$.

Thus $\rho\mapsto d_A(\rho')$ defines a homomorphism from $Z^1\(G,H^2(T;\Z)\)$
into the Moore group
$H^3\(G,C(T,\T)\)$.  If $[\rho]=0$ in $H^1$, then there is an
element $b\in H^2(T;\Z)$ such that $\rho_s =
b(s\cdot b)^{-1}$.  We can lift $b$ to an element $\theta
\in \Autc\(C_0(T,\K)\)$, and define
$\alpha_s=\theta\ttau_s\theta^{-1}$.  Then
\begin{equation*}
\alpha_s\alpha_t=
(\theta\ttau_s\theta^{-1})(\theta\ttau_t\theta^{-1})=
\theta\ttau_{st}\theta^{-1}
=\alpha_{st}.
\end{equation*}
Thus $\rho'$ lifts to a (trivially) twisted action, our homomorphism
factors through $H^1\(G,H^2(T;\Z)\)$, and we can define the required
homomorphism
$d_2'$ by $d_2'\([\rho]\)=d_A(\rho')$.

Now suppose that $\c\in H^3(T;\Z)^G$ is in the kernel of $d_2$.
If $A$ is a stable continuous-trace \cs-algebra with spectrum $T$
and with $\delta(A)=\c$, then $d_2(\c)=0$ implies that
there is a homomorphism $\g:G\to\Out(A)$ lifting the canonical
map $\ell:G\to\Homeo_\c(T)$ (i.e., $h\circ \g=\ell$).
Then \lemref{lem:5.5}(1) gives us an obstruction $d_A(\g)$.  If
$\delta:G\to\Out(A)$ also satisfies
$h\circ\delta=\ell$, then we claim that $d_A(\g)d_A(\delta)^{-1}$
belongs to $\Im(d_2')$.  To see this, first recall that, as
pointed out in the beginning of this section, there are
$C_0(T)$-isomorphisms $\phi_1:A\Ttensor \bA \to
A\Ttensor \bA\tensor\K$ and $\phi_2:A\Ttensor \bA\tensor\K
\to C_0(T,\K)$.  (We have already seen that
$\delta(A\Ttensor\bA)=0$.)  Thus,
\begin{align*}
d_A(\g)d_A(\delta)^{-1} &= d_A(\g)d_{\bA}(\bar \delta) \\
&= d_{A\Ttensor \bA}(\g\tensor\bar\delta )\\
&= d_{A\Ttensor \bA\tensor\K}\(\Ad(\phi_1)[\g\tensor\bar\delta]\)
\\
&=d_{C_0(T,\K)}\(\Ad(\phi_2\phi_1)[\g\tensor\bar\delta]\),
\end{align*}
which is by definition $d_2'(\rho)$ where $\rho_s=\zeta
\(\Ad(\phi_1\phi_2)[\gamma\tensor\bar\delta] \circ \tau_s^{-1}\)$.
This establishes the claim. Consequently, we may make the following
definition of~$d_3$.

\begin{definition}
Given $\c\in\ker d_2$, then $d_3(\c)$ is defined to be the
class of $d_A(\g)$ in $H^3\(G,C(T,\T)\)$ modulo the image of
$d'_2$, where $A$ is a stable continuous-trace \cs-algebra with
spectrum $T$ such that $\delta(A)=\c$, and $\g$ is any lift of the
canonical map $\ell:G\to\Homeo_\c(T)$ to a homomorphism $\g:G
\to\Out(A)$.
\end{definition}
Notice that it follows from \lemref{lem:5.5}(2) and
\propref{prop:2.2} that $d_3$ is a homomorphism.
Now we are ready to identify the kernel of $d_3$ with $\Im(F)$.
\begin{lem}\label{lem:31}
Suppose that $\c\in H^3(T;\Z)^G$.  Then $\c$
is in the kernel of both $d_2$ and
$d_3$ if and only if there there is a twisted action $\alpha:
G\to \Aut(A)$, $u:G\times G\to\UM(A)$ such that $h(\alpha_s)=\ell_s$
for all $s\in G$.
\end{lem}
\begin{pf}
We prove only the non-trivial direction.
Since $d_2(\c)=0$, there is a homomorphism $\g:G\to\Out(A)$ such
that $h\circ\g=\ell$, where $A$ is a stable, continuous-trace
\cs-algebra with spectrum $T$ and $\delta(A)=\c$.
Since $d_3(\c)=0$, $d_A(\g)\in\Im(d_2')$.
Thus we can choose $\rho\in Z^1\(G,H^2(T;\Z)\)$ with $d_2'(\rho)=
d_A(\g)^{-1}$.  Since $A$ is stable, $A$ is $C_0(T)$-isomorphic
to $A\tensor\K$ by \cite[Lemma~4.3]{pr2}, and there is a
$C_0(T)$-isomorphism $\phi:A\Ttensor C_0(T,\K)\to A$.
As above, let
$\rho'_s=\zeta^{-1}(\rho_s)\circ\tau\in
\Out\(C_0(T,\K)\)$. Then
\begin{align*}
d_A\(\Ad(\phi)(\g\tensor\rho')\) &=
d_{A\Ttensor C_0(T,\K)}\(\g\tensor\rho') \\
&=d_A(\g)d_{C_0(T,\K)}(\rho')\qquad\text{by \lemref{lem:5.5}(2)}\\
\intertext{which by definition of $d_2'$ is}
&= d_A(\g)d_2'(\rho),
\end{align*}
which is trivial by construction.  The result now follows from
\lemref{lem:5.5}(1).
\end{pf}
\begin{thm} \label{thm:5.8}
Suppose
that $\ell:G\to\homeot$ is induced by a
second countable locally compact
transformation group $(G,T)$ with $H^2(T;\Z)$ countable.
Then the image of the Forgetful Homomorphism $F$ is exactly those
classes $\c\in H^3(T;\Z)$ which lie in $H^3(T;\Z)^G$ and which
satisfy $d_2(\c)=0=d_3(\c)$.  In particular, a stable
continuous-trace \cs-algebra $A$ with spectrum $T$ and
Dixmier-Douady class $\delta(A)$
admits an action $\alpha:G\to\Aut(A)$ covering $\ell$
if and only if $\delta(A)\in H^3(T;\Z)^G$, $d_2\(\delta(A)\)=0$,
and $d_3\(\delta(A)\)=0$.
\end{thm}
\begin{pf}
The first statement follows from the second.
Furthermore, the
necessity of these conditions is evident.
So suppose that $d_3\(\delta(A)\)=0$.
By \lemref{lem:31} there is a twisted action $(\beta,u)$ on $A$
such that $h(\beta_s)=\ell_s$.  Using the stabilization trick
\cite[Theorem~3.4]{para1}, we note that $\beta\tensor i$ is exterior
equivalent (see \cite[Definition~3.1]{para1}) to an (ordinary)
action $\alpha$ on $A\tensor\K$.  Since $A\tensor\K$ is
$C_0(T)$-isomorphic to $A$, we are done once we show that
$h(\alpha_s)=h(\beta_s)$.  But, for each fixed $s$, $\beta_s\tensor
i$ differs from $\alpha_s$ by an inner automorphism $\Ad v_s$.
Then if $\pi\in \hA$,
\begin{align*}
s^{-1}\cdot(\pi\tensor i)
&= \pi\tensor i\circ \alpha_s =
\pi\tensor i\circ\Ad v_s \circ \beta_s\tensor i \\
&=\pi\tensor i(v_s)\(\pi\tensor i\circ \beta_s \tensor
i\)\pi\tensor i(v_s^*) \\
&\sim \pi\tensor i \circ \beta_s\tensor i =(s^{-1} \cdot \pi)
\tensor i,
\end{align*}
which completes the proof.
\end{pf}

%
%   SCCS Info: ckrw6.tex version 3.1: 10/19/93
%
\section{The Structure Theorem}
\label{sec:structure}

\begin{thm}
\label{thm:structure}
Suppose $(G,T)$ is a second countable locally compact transformation
group such that $H^2(T;\Z)$ is countable. There are homomorphisms
\begin{align*}
d_2&:H^3(T;\Z)^G\to H^2\(G,H^2(T;\Z)\),\\
d_2'&:H^1\(G,H^2(T;\Z)\)\to H^3\(G,C(T,\T)\),\\
d_2''&:H^2(T;\Z)^G\to
H^2\(G,C(T,\T)\),\text{ and}\\
d_3&:\ker(d_2)\to H^3\(G,C(T,\T)\)/\Im(d_2'),
\end{align*}
with the following properties.  (Indeed, $d_2$, $d_2'$, and $d_3$
are the homomorphisms defined in the previous section.)
\begin{enumerate}
\item The homomorphism $F:[A,\alpha]\mapsto\delta(A)$ of $\brgt$ into
$H^3(T,\Z)^G$ has range $\ker(d_3)$, and kernel consisting of all
classes of the form $[C_0(T,\K),\alpha]$.
\item Let $\zeta:\Autc{C_0(T,\K)}\to H^2(T;\Z)$ be the homomorphism
of \cite[Theorem~2.1]{pr1}. Then the homomorphism
$\eta:\ker(F)\to  H^1\(G,H^2(T;\Z)\)$ defined by
$$
\eta\big(C_0(T,\K),\alpha\big)(s)=\zeta(\alpha_s\circ\tau^{-1}_s),
$$
has range $\ker(d_2')$.
\item
For each cocycle $\w\in Z^2\(G,C(T,\T)\)$, choose a Borel map $u:G\to
UM\(C_0(T,\K)\)$ satisfying
\begin{equation}\label{eq:oldhat}
u_s\tau_s(u_t)=\w(s,t)u_{st},
\end{equation}
and define $\xi(\w)=[C_0(T,\K),\Ad{u}\circ\tau]$. Then $\xi$ is a
well-defined homomorphism of $H^2\(G,C(T,\T)\)$ into $\brgt$ with
$\Im(\xi)=\ker(\eta)$ and $\ker(\xi)=\Im(d_2'')$.

\end{enumerate}

\end{thm}

\begin{pf}
It follows from \thmref{thm:5.8} that the image of $F$ is $\ker(d_3)$.
If $(A,\alpha)\in\ker(F)$, then $\delta(A)=0$ and $A\tensor\K$ is
$C_0(T)$-isomorphic to $C_0(T,\K)$.  Thus
\begin{equation*}
(A,\alpha)\sim(A\tensor\K,\alpha \tensor i)\sim\(C_0(T,\K),\alpha'\),
\end{equation*}
where $\alpha'=\theta\circ(\alpha\tensor i)\circ\theta^{-1}$ for
some $C_0(T)$-isomorphism $\theta$.
This proves part~(1).

If $[A,\alpha]=[C_0(T,\K),\alpha']\in\ker(F)$,
then we
want to define
\begin{equation}\label{eq:defeta}
\eta\([A,\alpha]\)= \bigl[
s\mapsto \zeta_{C_0(T,\K)}(\alpha_s'\circ\ttau_s^{-1})\bigr],
\end{equation}
and, we need to verify that $\eta$ is well defined.  So, suppose
that $(A,\alpha)\sim\(C_0(T,\K),\beta\)$ as well;  $\alpha'$
and $\beta$ are outer conjugate over $T$ by \lemref{lem:combes}.
But if $\alpha'$ is exterior
equivalent to $\beta$, say $\beta_s=\Ad(u_s)\circ \alpha_s'$, then
\begin{equation*}
\zeta\(\(\Ad(u_s)\circ \alpha_s'\)\circ\ttau_s^{-1}\)=
\zeta\(\Ad(u_s)\circ(\alpha_s'\circ \ttau_s^{-1})\)=
\zeta(\alpha_s'\circ\ttau_s^{-1}).
\end{equation*}
On the other
hand, if $\beta$ is conjugate to $\alpha'$, say $\beta_s=\Phi\circ
\alpha_s'\circ\Phi^{-1}$ with $\Phi\in\Autc\(C_0(T,\K)\)$, then
\begin{align}
\zeta(\Phi\circ\alpha_s'\circ\Phi^{-1}\circ\ttau_s^{-1})
&=\zeta(\Phi\circ\alpha_s'\circ\ttau_s^{-1}\circ \ttau_s\circ
\Phi^{-1}\circ\ttau_s^{-1})\notag\\
\intertext{which, since the range of $\zeta$ is abelian, is}
&=
\zeta\((\alpha_s'\circ\ttau_s^{-1})\circ \ttau_s
\circ\Phi^{-1}\circ\ttau_s^{-1}\circ\Phi\)\notag\\
&=
\zeta(\alpha_s'\circ\ttau_s^{-1}) \zeta(\ttau_s
\circ\Phi^{-1}\circ\ttau_s^{-1}) \zeta(\Phi) \notag \\
\intertext{which, since $\zeta$ is equivariant, is}
&=
\zeta(\alpha_s'\circ\ttau_s^{-1})
s\cdot\zeta(\Phi)^{-1}\zeta(\Phi) .\label{eq:last}
\end{align}
Thus image of $\alpha$ and $\beta$ differ by a coboundary in $B^1\(
G,H^2(T;\Z)\)$ and \eqref{eq:defeta} gives a well-defined map of
$\ker(F)$ into $H^1\(G,H^2(T;\Z)\)$.
It is not hard to check, using \cite[Proposition~3.10]{pr2},
that $\eta$ is a homomorphism.
Furthermore, $d_2'\([\rho]\)=0$ if and only if there is a twisted
action $(\alpha,u)$ on $C(X,\K)$ such that
$\zeta(\alpha_s\circ\ttau_s^{-1})=\rho_s$.  Thus if $[\rho]=\eta(A,
\alpha)$, then certainly $d_2'\([\rho]\)=0$.  On the other hand,
if $d_2'\([\rho]\)=0$, then let $(\alpha,u)$ be a twisted action
with $\zeta(\alpha\circ\ttau^{-1})=\rho$.  Then the stabilization
trick \cite[Theorem~3.4]{para1} implies that there is an action
$\beta$ of $G$ on $C_0(T,\K)\tensor\K$ which is exterior
equivalent to $\alpha\tensor i$: say $\beta_s=\Ad(v_s)\circ
(\alpha_s\tensor i)$.
Then $\beta_s\circ\ttau_s^{-1}=\Ad(v_s)\circ (\alpha_s\tensor
i)\circ \ttau_s^{-1}=\Ad(v_s)\circ(\alpha_s\circ\ttau_s^{-1})\tensor
i$, so
\begin{equation*}
\zeta(\beta_s\circ\ttau_s^{-1})=
\zeta(\alpha_s\circ\ttau_s^{-1}\tensor
i)=\zeta(\alpha_s\circ\ttau_s^{-1}) = \rho_s.
\end{equation*}
Since $C_0(T,\K)\tensor\K$ is $C_0(T)$-isomorphic to $C_0(T,\K)$,
say by $\phi$, we have
$\eta\(C_0(T,\K),\phi\circ\beta\circ\phi^{-1}\)=[\rho]$.  Thus,
the image of $\eta$ is equal to $\ker(d_2')$ as required.
This proves~(2).

For convenience, let $A=C_0(T,\K)$.
To define $\xi$, we first need to note
that for every $\w\in Z^2\(G,C(T,\T)\)$ there is a Borel map $u:G
\to\UM(A)$ satisfying \eqref{eq:oldhat}.
However, it was shown in the proof of \cite[Proposition~3.1]{horr} that
\begin{equation*}
\(u^\w_t(x)f\)(s)=\w(s,t)(s\cdot x)f(st)\qquad\text{for $f\in L^2(G)$}
\end{equation*}
gives for each $t\in G$ a strongly continuous map $u^\w_t:
T\to\U\(L^2(G)\)$, defining a unitary multiplier $u^\w_t$ of $
C_0\(T,\K\(L^2(G)\)\)$, and that
$u^\w:G\to\UM\(C_0\(T,\K\(L^2(G)\)\)\)$ is then a Borel map
satisfying \eqref{eq:oldhat}.
If $G$ is infinite, then  $u=u^\w$ is the desired map.  Otherwise, we can
stabilize and let $u=u^\w\tensor i$.
In any case, it is clear that $(A,\Ad(u)\circ\tau)$ is an element of
$\sbrgt$.  If $v:G\to\UM(A)$ also satisfies \eqref{eq:oldhat}, then
$s\mapsto u_sv^*_s$ gives an exterior equivalence between
$\Ad(u)\circ\tau$ and $\Ad(v)\circ\tau$, so $[A,\Ad(u)\circ\tau]$ is
independent of the choice of $u$.
Since we can absorb a coboundary in $B^2\(G,C(T,\T)\)$ into the unitary
$u$ without changing $\Ad(u)\circ\tau$, we have a well-defined class
$\xi\([\w]\)$ in $\brgt$, and another routine argument shows that $\xi$
is a homomorphism.

To see that the image of $\xi$ is the kernel of $\eta$, note that
$(\Ad(u)\circ\tau)\circ\tau^{-1}$ consists of inner automorphisms, and
hence $\eta\circ\xi\([\w]\)=\eta\(A,\Ad(u)\circ\tau\)$ is identically
zero.
On the other hand, if $\eta(A,\alpha)=0$, then there exists
$\phi\in\Autc(A)$ such that $\zeta(\alpha_s\circ\tau_s^{-1})=\zeta(\phi)
s\cdot\zeta(\phi)^{-1}=\zeta(\tau_s\circ\phi^{-1}\circ\tau_s^{-1})$ for
all $s\in G$.  Thus $\phi^{-1}\circ\alpha_s\circ\phi$ differs from
$\tau_s$ by inner automorphisms $\Ad(u_s)$; we can choose $u:G\to\UM(A)$
to be Borel, $u$ satisfies \eqref{eq:oldhat} for some cocycle $\w\in
Z^2\(G,C(T,\T)\)$, and $\phi^{-1}$ is an isomorphism taking $(A,\alpha)$
to $(A,\Ad(u)\circ\tau)$.  Thus,
$[A,\alpha]=[A,\Ad(u)\circ\tau]=\xi\([\w]\)$ lies in the image of $\xi$.

To define the homomorphism $d_s''$, we choose an automorphism
$\phi\in\Autc(A)$ such that $\zeta(\phi)\in H^2(T;\Z)^G$, which means
precisely that $[\phi]=s\cdot[\phi]:=[\tau_s\circ\phi\circ\tau_s^{-1}]$.
Then $s\mapsto\tau_s\circ\phi\circ\tau_s^{-1}\circ\phi^{-1}$ is a Borel
map of $G$ into $\Inn(A)$, and there is a Borel map $u:G\to\UM(A)$ such
that $\Ad(u_s)\circ\tau_s\circ\phi=\phi\circ\tau_s$.  The usual argument
shows that $\Ad\(u_s\tau_s(u_t)\)=\Ad(u_{st})$, so $u$ satisfies
\eqref{eq:oldhat} for some cocycle $\w\in Z^2\(G,C(T,\T)\)$, and $\phi$
gives an isomorphism of $(A,\tau)$ onto the system
$\(A,\Ad(u)\circ\tau\)$ representing $\xi\([\w]\)$.
The choice of $u$ was unique up to multiplication by a function
$\rho:G\to C(T,\T)$, so $\w$ is unique up to multiplication by
$\partial\rho$; choosing a different representative $\Ad(v)\circ\phi$ for
$\zeta(\phi)$ would force us to use $v^*u_s\tau_s(v)$ in place of $u_s$,
which would not change $\w$.  Thus we have a well-defined class
$[\w]=d_2''\(\zeta(\phi)\)$ in $H^2\(G,C(T,\T)\)$.  If
$\Ad(u_s)\circ\tau_s\circ\phi=\phi\circ\tau_s$ and
$\Ad(v_s)\circ\tau_s\circ\psi=\psi\circ\tau_s$, then
\begin{equation*}
\Ad\(\phi(v_s)u_s\)\circ\tau_s\circ(\phi\circ\psi)=
(\phi\circ\psi)\circ\tau_s,
\end{equation*}
and a routine calculation using this shows that $d_2''\(\zeta(\phi\circ
\psi)\)=d_2''\(\zeta(\phi)\)+d_2''\(\zeta(\psi)\)$.

We have already seen that
\begin{equation*}
\xi\(d_2''\(\zeta(\phi)\)\)=[A,\Ad(u)\circ\tau]=[A,\tau],
\end{equation*}
so $\xi\circ d_2''=0$.  Conversely, if $\xi\([\w]\)=[A,\Ad(u)\circ\tau]=0$
in $\brgt$, then \lemref{lem:combes} gives a $C_0(T)$-automorphism of $A$
such that $\phi^{-1}\circ\(\Ad(u)\circ\tau)\circ\phi$ is exterior
equivalent to $\tau$.  But if $v$ is a $\tau$-1-cocycle such that
\begin{equation}\label{eq:idag}
\phi^{-1}\circ\(\Ad(u_s)\circ\tau_s)\circ\phi=\Ad(v_s)\circ\tau_s,
\end{equation}
then $\Ad\(\phi(v_s^*)u_s\)\circ\tau_s\circ\phi=\phi\circ\tau_s$ and a
quick calculation using \eqref{eq:idag} shows that
$$
\(\phi(v_s^*)u_s\)\tau_s\(\phi(v_t^*)u_t\)=\w(s,t)\phi(v_{st}^*)u_{st},
$$
so that $[\w]=d_2''\(\zeta(\phi)\)$.
\end{pf}

%
%   SCCS Info: ckrw7.tex version 3.1: 10/19/93
%
\section{Examples and Applications}
\label{sec:7}

\subsection{Actions of $\R$.}
When $G=\R$, we can sharpen the conclusion of \thmref{thm:5.8}
considerably, and we obtain the generalization of
\cite[Theorem~4.12]{rr} mentioned in the introduction.

\begin{cor}\label{cor:5.9}
Suppose that $(\R,T)$ is a second countable locally compact
transformation group with $H^1(T;\Z)$ and
$H^2(T;\Z)$ countable, and $A$ is a stable
continuous-trace \cs-algebra with spectrum $T$. Then there is
always, up to exterior
equivalence, exactly one action $\alpha:\R\to\Aut(A)$ covering
the given action on $T$.
\end{cor}

\begin{pf} Since there is an action on $A$ covering the given action on $T$
\iff{} $\delta(A)$ belongs to the range of $F$, we have to prove that
$F$ is surjective. By \thmref{thm:5.8}, this is equivalent to proving
that $H^3(T;\Z)^{\R}=H^3(T;\Z)$, that $d_2=0$, and that $d_3=0$. The
connectedness of $\R$ implies that $l_s$ is homotopic to
$l_e=\id$, and consequently that $(l_s)^*=\id$ for all $s$,
giving $H^3(T;\Z)^{\R}=H^3(T;\Z)$. The homomorphism $d_2$ is 0 because
$H^2(\R,M)$ is trivial for any discrete $\R$-module $M$ \cite[Theorem
4]{wigner}, and we are assuming that $M=H^2(T,\Z)$ is countable.
Finally, Theorem 4.1 of \cite{rr} says that $H^3(\R,C(T,\T))=0$, and
therefore $d_3$ is also 0. \end{pf}

\subsection{Free and Proper Actions.}\label{sec:7.2}
 We now suppose
that $G$ acts freely and properly on $T$; in other words, that $s\cdot
x=x$ \iff{} $ s=e$, and that $(s,x)\mapsto(s\cdot x,x)$ is a
proper map of $G\times T$ into $T\times T$. (If $G$ is a Lie group,
it follows from a theorem of Palais \cite[Theorem~4.1]{palais} that
these are precisely the locally trivial principal $G$-bundles.) Let
$p:T\to T/G$ denote the orbit map. By \cite[Theorem~1.1]{rr}, every
$(A,\alpha)\in\sbrgt$ in which $A$ is stable is equivalent to
an element of the form
$(p^*B,p^*\id)=\big(C_0(T)\otimes_{C(T/G)}B,\tau\otimes_{C_0(T/G)}\id
\big)$---indeed,
we can take $B$ to be the crossed product $A\rtimes_\alpha G$. Thus the
orbit map $p$ induces a homomorphism $p^*:B\mapsto (p^*B,p^*{\id})$ of
$\br(T/G)$ onto $\brgt$. Since the crossed product map
$\rtimes_\alpha:(A,\alpha)\mapsto [A\rtimes_\alpha G]$ is well-defined
on $\brgt$ by the Combes--Curto-Muhly-Williams Theorem, and
$$
p^*B\rtimes_{p^*{\id}}G\cong\big(C_0(T)\rtimes_\tau
G\big)\otimes_{C(T/G)}B\cong C_0(T/G,\K)\otimes_{C(T/G)}B\cong
B\otimes\K,
$$
the map $\rtimes_\alpha$ is an inverse for $p^*$, and $p^*$ is an
isomorphism.

\subsection{Trivial Actions.} If $G$ acts trivially on $T$, then
$F:\brgt\to \br(T)$ is trivially surjective (given $A$, take
$(A,\id)$), and a quick look at the definition in
\secref{sec:structure}
 shows that
the map $d_2'':H^2(T;\Z)^G\to H^2(G,C(T,\T))$ is zero. Thus our
structure theorem gives an exact sequence
\begin{ediagram*}
\node{0} \arrow{e} \node{H^2(G,C(T,\T))} \arrow[2]{e,t}{\xi}
\node[2]{\ker(F)=\set{\big[C_0(T,\K),\alpha\big]}}
\arrow{wsw,l}{\eta} \\
\node[2]{\operatorname{Hom}\(G,H^2(T,\Z)\)}\arrow[2]{e,t}{d_2'}
\node[2]{ H^3(G,C(T,\T)).}
\end{ediagram*}
If $G$ is also connected, then every action
$\alpha:G\to\Autc{C(T,\K)}$ has range lying in the open subgroup
$\Inn{C(T,\K)}$, so $\eta=0$ and $H^2(G,C(T,\T))$ classifies the
actions of $G$ on $C(T,\K)$ (and, by Remark~\ref{rem:4.10}, any other
stable continuous-trace algebra with spectrum $T$\/); see
\cite[\S0]{rr}. On the other hand, if $G=\Z$, $H^2(G,C(T,\T))=0$, and
$\eta$ is an isomorphism of $\ker(F)\cong\Autc{C_0(T,\K)}$ onto
$\operatorname{Hom}(G,H^2(T;\Z))=H^2(T;\Z)$ by \cite[Theorem
2.1]{pr1}; if $G=\R$, all the groups in the sequence vanish, and
$\ker(F)$ is trivial.

There are two extreme special cases.  When $G$ is trivial, we
recover the Dixmier-Douady isomorphism $\brt\cong H^3(T;\Z)$ (see
Remark~\ref{rem:3.3}).  When the space $T$ consists of a single
point, the elements of $\sbrgt$ are systems $(\K,G,\alpha)$, and
we recover from \thmref{thm:structure} the parameterization of
actions of $G$ on $\K$ by the Moore cohomology group $H^2(G,\T)$.

\subsection{$N$-proper systems.} Suppose $G$ is abelian and there is
a closed subgroup $N$ such that $s\cdot x=x$ \iff{} $x\in N$, so
that $G/N$ acts freely on $T$. We suppose that $p:T\to T/G$ is a locally
trivial principal $G/N$-bundle. If $N$ is a compactly generated group
with $H^2(N,\T)=0$, and $(A,\alpha)\in\sbrgt$, then
$\alpha\restr N$ is locally unitary
\cite[Corollary~2.2]{ros2}, so the system
$(A,\alpha)$ is $N$-principal in the sense of \cite{90a,90b,90c}.
Provided the quotient map $\widehat G\to\widehat N$ has local
cross-sections (equivalently, $\widehat G$ is locally trivial as an
$N^\perp$-bundle), then Proposition 4.8, Corollary 4.4 and Theorem 6.3
of \cite{90c} imply that the Dixmier-Douady invariant of \cite{90c}
induces an isomorphism of $\brgt$ onto the equivariant sheaf
cohomology group $H^2_G(T,\sT)$ studied in \cite{95a}. The Gysin
sequence of \cite{95a} then implies that we have an exact sequence
\begin{ediagram*}\dgARROWLENGTH=.2em
\arrow[2]{e}\node[2]{\br(T/G)}\arrow[4]{e,t}{p^*}
\node[4]{\brgt}\arrow[5]{e,t}{b}\node[5]{H^1(T/G,\widehat{\cal
N})}\arrow[5]{e} \node[5]{H^3(T/G,\sT)}\arrow[3]{e} \node[3]{.}
\end{ediagram*}
The homomorphism $p^*$ takes a continuous-trace algebra $B$ with
spectrum $T/G$ to $(p^*B,p^*\id)$, and the homomorphism $b$ takes
$(A,\alpha)$ to the class of the principal $\widehat N$-bundle
$\spec{A\rtimes_\alpha G}\to T/G$.  Thus taking $N=\{e\}$ gives a
special case of \secref{sec:7.2} concerning free actions, and taking
$G=N$ gives a short exact sequence
\begin{ediagram*}
\node{0}\arrow{e} \node{\br(T)}\arrow{e}
\node{\brgt}\arrow{e,t}{\zeta}
\node{H^1(T,\widehat G)} \arrow{e} \node{0}
\end{ediagram*}
summarizing the results of \cite{pr2} for the case where $G$ acts
trivially on $T$. Various other special cases are considered in the
last section of \cite{95a}. However, we stress that the group
$H^2_G(T,\sT)$ is {\em not} a true equivariant cohomology group in
the sense of \cite{groth}: if $H^2(N,\T)\not=0$, then the systems
$(A,\alpha)$ which are locally Morita equivalent to
$\big(C_0(T),\tau\big)$, and hence classifiable by their
Dixmier-Douady class in $H^2_G(T,\sT)$, form a proper subgroup of
$\brgt$.

\def\mathcs{C^{\displaystyle *}} \def\cs{\ifmmode\mathcs\else$\mathcs$\fi}
\ifx\undefined\bysame
\newcommand{\bysame}{\leavevmode\hbox to3em{\hrulefill}\,}
\fi

\end{document}